\documentclass[11pt]{article}
\usepackage[T1]{fontenc}
\usepackage[utf8]{inputenc}
\usepackage{amsmath,amssymb,amsthm,microtype}
\usepackage{array}
\usepackage[colorlinks=true,linkcolor=blue,citecolor=blue,urlcolor=blue]{hyperref}
\newtheorem{theorem}{Theorem}[section]
\newtheorem{lemma}[theorem]{Lemma}
\newtheorem{proposition}[theorem]{Proposition}
\newtheorem{corollary}[theorem]{Corollary}
\theoremstyle{definition}
\newtheorem{definition}[theorem]{Definition}
\newtheorem{remark}[theorem]{Remark}
\newtheorem{question}[theorem]{Question}
\newcommand{\R}{\mathbb R}
\newcommand{\Q}{\mathbb Q}
\newcommand{\Z}{\mathbb Z}
\newcommand{\T}{\mathbb T}
\newcommand{\nrm}[1]{\lVert #1\rVert}
\title{Joint Approximation and the Mixed-Prime Character Group}
\author{Irina Yaroshevskaya}
\date{}
\hypersetup{pdftitle={Joint Approximation and the Mixed-Prime Character Group},pdfauthor={Irina Yaroshevskaya}}
\begin{document}
\maketitle
\begin{abstract}
For the integer metric defined by the square denominators of two distinct primes, we prove simultaneous approximation to zero in one prime metric and to any prescribed integer in the other. A normalized convolution argument gives a uniform central spectral bound. Random relative translations of smooth weights produce an exact positive cascade; Fourier positivity places every normalized finer convolution below its value at zero. A far-field recursion, the growth of merged square powers, and a two-term logarithmic-form estimate give summable factor bounds. Consequently, for every fixed pair of distinct primes $p,q$, the mixed continuous character group is $\Z[1/(pq)]/\Z$: no irrational parameter is continuous. The mixed completion is the product of the one-prime completions. For each fixed pair this also gives a compatible rational join and a strictly finer compatible join of the associated real TVG topologies, with the same prescribed real character group. Exact witness search terminates at every strict precision and integer target. We retain the structural, density, divisibility, sequence and finite-certificate results leading to this proof. The existence argument does not give practical witness-size or runtime bounds, and it does not assume independence or equidistribution of the two digit systems. The joint approximation theorem has been formalized in Lean and Mathlib conditionally on the two-term logarithmic separation estimate \eqref{eq:two-term-separation}, which is an explicit hypothesis of the formal statement; that estimate, a consequence of Matveev's bound for linear forms in logarithms, is not formalized. The corresponding problem for three or more primes remains separate.
\end{abstract}

\clearpage
\tableofcontents
\clearpage

\section{Relation to the earlier articles}
The articles \cite{Previous,HamelPair} study topologies on the real line generated by a rational-valued Hamel function. Section~7 of Part~I of \cite{HamelPair} isolates a problem about the group of integers $\Z$, for the primes $2$ and $3$. Write $\T=\R/\Z$ and let $\nrm{t}$ denote the distance from a real number $t$ to the nearest integer, so $0\le\nrm{t}\le 1/2$. For a finite nonempty set $S$ of primes and an integer $k$ put
\[
P_S(k)=\max_{\ell\in S}\ \sup_{n\ge1}\nrm{k/\ell^{\,n^2}}.
\]
Part~I of \cite{HamelPair} uses $S=\{2,3\}$. Throughout this article $S$ is a fixed finite nonempty set of primes, and we write $P=P_S$; the case $S=\{2,3\}$ is the one of Part~I of \cite{HamelPair}. Let $s$ be the product of the primes in $S$, and let $\Z[1/S]$ be the subring of $\Q$ generated by the fractions $1/\ell$ with $\ell\in S$, so that $\Z[1/S]=\Z[1/6]$ for $S=\{2,3\}$.

A \emph{character} of $\Z$ is a homomorphism $\Z\to\T$; every character has the form $k\mapsto k\alpha+\Z$ for a unique $\alpha\in\T$, which we call its \emph{parameter}. Let $H$ be the set of parameters $\alpha$ whose character is continuous for the metric $d(k,k')=P(k-k')$ on $\Z$. For $S=\{2,3\}$, Theorem~7.1 of Part~I of \cite{HamelPair} shows that $H$ has Hausdorff dimension zero and contains $\Z[1/6]/\Z$, Theorem~\ref{thm:joint-approximation} and Corollary~\ref{cor:JA-topology} below determine this group exactly; Corollary~\ref{rem:JA-other-pairs} treats every fixed pair of distinct primes. For general $S$ the inclusion $\Z[1/S]/\Z\subseteq H$ is proved in Remark~\ref{rem:inclusion}.

The principal new argument is in Section~\ref{sec:normalized-product}: an exact normalized product establishes a scale-uniform spectral bound. Section~\ref{sec:random-JA} constructs the positive cascade by random relative translations, controls its spatial masses, and transfers that bound to joint approximation. Section~\ref{sec:JA-consequences} determines the character group and completion and records the approximation consequences. The joint approximation theorem has been formalized in Lean and Mathlib conditionally on the two-term separation estimate \eqref{eq:two-term-separation}; the scope is specified in the section on checks by Lean and Mathlib below.

The preceding sections supply the structural background: continuity and flatness, rational open subgroups, null sequences, one-prime duals, counting and stationary digit models, splitting under joint approximation, certified interval search, quotient compactness, density criteria, and collision-energy comparisons. The proof uses a logarithmic-form theorem of Matveev as well as finite Fourier positivity. Character separation for the quotient criterion uses Pontryagin duality. The earlier transversality results remain supplementary and are not used to infer fine-scale independence. Section~\ref{sec:open} lists the remaining quantitative and finite-family questions.

The individual square-power topologies are instances of the divisibility-chain construction of Au\ss enhofer and de la Barrera Mayoral \cite[Theorem~4.4]{ABM}. Applying that theorem to the single chain $b_n=(pq)^{n^2}$ also gives a topology with dual $\Z[1/(pq)]/\Z$. The arithmetic issue here is the dual of the specified join of the two separate prime-square topologies, together with joint approximation and the resulting completion. General criteria for compatible joins and extension of characters from a diagonal subgroup belong to the established theory; see \cite[Proposition~3.9 and Theorem~3.10]{GabCompatible}. Related compatible suprema of sparse subsequence topologies were studied in \cite[Section~3]{DBM2017}. These general constructions are distinct from the mixed-prime approximation problem considered here.

\section{The pseudo-norm $P$}\label{sec:pnorm}
For $\delta>0$ put $B_\delta=\{k\in\Z:P(k)<\delta\}$.

\begin{lemma}\label{lem:basic}
For all integers $k,k'$ and $j\ge1$:
\begin{enumerate}
\item $P(0)=0$, $P(-k)=P(k)$, and $P(k)>0$ if $k\ne0$;
\item $P(k+k')\le P(k)+P(k')$ and $P(jk)\le jP(k)$;
\item $jB_{\delta/j}\subseteq B_\delta$;
\item if $\delta\le \ell^{-j^2}$ for every $\ell\in S$, then $B_\delta\subseteq s^{j^2}\Z$; if $\delta\le1/\ell$ for every $\ell\in S$, then $B_\delta\subseteq s\Z$.
\end{enumerate}
\end{lemma}
\begin{proof}
(1) Evenness and $P(0)=0$ are clear. If $k\ne0$ and $\ell\in S$, then for $n$ with $\ell^{n^2}>2|k|$ we have $\nrm{k/\ell^{n^2}}=|k|/\ell^{n^2}>0$.

(2) The map $t\mapsto\nrm{t}$ satisfies $\nrm{t+t'}\le\nrm{t}+\nrm{t'}$, hence each $k\mapsto\nrm{k/\ell^{n^2}}$ is subadditive, and so are its suprema and their maximum $P$. The bound for $jk$ follows by induction.

(3) is a restatement of $P(jk)\le jP(k)$.

(4) Let $\ell\in S$ and take $n=j$. If $\nrm{k/\ell^{j^2}}<\delta\le\ell^{-j^2}$, then, because the numbers $k/\ell^{j^2}$ lie in $\ell^{-j^2}\Z$, the distance from $k/\ell^{j^2}$ to the nearest integer is a multiple of $\ell^{-j^2}$ and is smaller than $\ell^{-j^2}$, hence zero. So $\ell^{j^2}\mid k$. Distinct primes are coprime, so $s^{j^2}\mid k$. For $n=1$: $\nrm{k/\ell}$ is either $0$ or at least $1/\ell$, so $\nrm{k/\ell}<\delta\le1/\ell$ forces $\ell\mid k$, and then $s\mid k$.
\end{proof}

By Lemma~\ref{lem:basic}(1)--(2), $d(k,k')=P(k-k')$ is a translation-invariant metric, and $\{B_\delta\}$ is a neighbourhood base of $0$ for a group topology on $\Z$. Call it the \emph{$P$-topology}. A sequence $(k_j)$ is \emph{null} if $P(k_j)\to0$.

\begin{lemma}\label{lem:multiples}
Let $\theta\in\R$ and $j\ge1$. If $\nrm{i\theta}\le1/4$ for every $i=1,\dots,j$, then $\nrm{\theta}\le 1/(4j)$.
\end{lemma}
\begin{proof}
Put $t=\nrm\theta$ and suppose $t>1/(4j)$. Then $\theta$ differs from $t$ or from $-t$ by an integer, so $\nrm{i\theta}=\nrm{it}$ for every integer $i$. If $t=1/2$ the hypothesis already fails for $i=1$, so $t<1/2$. Let $i=\lfloor 1/(4t)\rfloor+1$. Then $i\le j$, because $1/(4t)<j$ and $j$ is an integer, and $1/4<it\le 1/4+t<3/4$. Hence $\nrm{it}=\min\{it,1-it\}>1/4$, a contradiction.
\end{proof}

\begin{proposition}[quarter-arc criterion]\label{prop:quarter}
Let $\alpha\in\T$. The character of $\alpha$ is continuous for the $P$-topology if and only if there is $\delta>0$ with $\nrm{k\alpha}\le1/4$ for all $k\in B_\delta$.
\end{proposition}
\begin{proof}
Necessity is the definition of continuity at $0$ with $\varepsilon=1/4$. For sufficiency let $j\ge2$ and $k\in B_{\delta/j}$. Then $ik\in B_\delta$ for $1\le i\le j$ by Lemma~\ref{lem:basic}(3), so $\nrm{ik\alpha}\le1/4$ for these $i$, and Lemma~\ref{lem:multiples} gives $\nrm{k\alpha}\le1/(4j)$. As $j$ is arbitrary, the character is continuous at $0$, hence continuous.
\end{proof}

\begin{lemma}\label{lem:nondiscrete}
For every finite set $S$ of primes and every $\delta>0$, $B_\delta$ contains a nonzero integer. Hence the $P$-topology is not discrete.
\end{lemma}
\begin{proof}
If $\delta>1/2$, every integer lies in $B_\delta$, because $P(k)\le1/2$. Let $\delta\le1/2$ and put $m=\lfloor1/\delta\rfloor+1$, so that $1/m<\delta$. For $K\ge1$ let $E_K$ be the set of pairs $(\ell,n)$ with $\ell\in S$, $n\ge1$ and $\ell^{n^2}\le K/\delta$, and let $d$ be its size; then $d\le|S|\sqrt{\log_2(K/\delta)}$. Send each integer $k\in\{0,1,\dots,K\}$ to the vector of the numbers $\lfloor m\{k/\ell^{n^2}\}\rfloor$, $(\ell,n)\in E_K$, where $\{t\}$ is the fractional part of $t$. This vector has at most $m^d$ values, and $m^d=\exp(O(\sqrt{\log K}))<K+1$ for large $K$. So two integers $0\le k_1<k_2\le K$ have the same vector. Put $k=k_2-k_1\in[1,K]$. For $(\ell,n)\in E_K$ the fractional parts of $k_1/\ell^{n^2}$ and $k_2/\ell^{n^2}$ lie in one interval of length $1/m$, hence $\nrm{k/\ell^{n^2}}<1/m<\delta$. For $\ell\in S$ and $n$ with $\ell^{n^2}>K/\delta$ we have $0<k/\ell^{n^2}<\delta$, so $\nrm{k/\ell^{n^2}}<\delta$, and these terms decrease as $n$ grows. Thus $P(k)<\delta$.
\end{proof}
\begin{proposition}[Lipschitz criterion]\label{prop:lipschitz}
Let $x\in\T$. The character of $x$ is continuous for the $P$-topology if and only if there is a constant $C>0$ with $\nrm{kx}\le C\,P(k)$ for all $k\in\Z$. Hence $H=\{x\in\T:\ \sup_{k\ne0}\nrm{kx}/P(k)<\infty\}$.
\end{proposition}
\begin{proof}
If $\nrm{kx}\le C\,P(k)$ for all $k$, then $P(k)<\varepsilon/C$ gives $\nrm{kx}<\varepsilon$, so the character is continuous. Conversely, by Proposition~\ref{prop:quarter} choose $\delta\in(0,1/2]$ with $\nrm{kx}\le1/4$ for all $k\in B_\delta$. Let $k\ne0$ and $p=P(k)>0$. If $p\ge\delta$, then $\nrm{kx}\le1/2\le p/(2\delta)$. If $p<\delta$, let $j\ge1$ be the largest integer with $jp<\delta$; then $j\ge\delta/(2p)$. By Lemma~\ref{lem:basic}(3), $ik\in B_\delta$ for $1\le i\le j$, so $\nrm{ikx}\le1/4$ for these $i$, and Lemma~\ref{lem:multiples} gives $\nrm{kx}\le1/(4j)\le p/(2\delta)$. Thus $C=1/(2\delta)$ works.
\end{proof}

\begin{theorem}[Flatness of continuous characters]\label{thm:flatness}
For every finite nonempty set $S$ of primes,
\[
x\in H
\quad\Longleftrightarrow\quad
\lim_{\delta\downarrow0}\ \sup_{0<P(k)<\delta}
\frac{\nrm{kx}}{P(k)}=0.
\]
Equivalently, continuity implies $\nrm{kx}=o(P(k))$ as nonzero integers $k$ tend to zero in the $P$-topology.
\end{theorem}
\begin{proof}
Let $x\in H$, and choose $C>0$ with $\nrm{kx}\le CP(k)$ by Proposition~\ref{prop:lipschitz}. Suppose the limit is not zero. There are $\kappa>0$ and a sequence of nonzero integers, null in $P$, such that $\nrm{k_jx}\ge\kappa P(k_j)$. Put $\eta=1/(16C)$. We select a subsequence $k_i$ recursively, together with integers $0=N_0<N_1<\cdots$.

After $N_{i-1}$ has been chosen, take $k_i$ sufficiently far along the witness sequence that, writing $p_i=P(k_i)$,
\[
p_i<\eta/2,\qquad
p_i<\min_{\ell\in S}\ell^{-N_{i-1}^2}.
\]
For $\ell\in S$ and $n\le N_{i-1}$, a nonzero residue of $k_i/\ell^{n^2}$ has norm at least $\ell^{-n^2}$. Hence all these early coordinates are exactly zero modulo one. Set
\[
t_i=\lfloor\eta/p_i\rfloor,\qquad v_i=t_i k_i.
\]
Then $t_i\ge\eta/(2p_i)$, $t_ip_i\le\eta$, and $P(v_i)\le\eta$. Once $v_i$ is fixed, choose $N_i>N_{i-1}$ such that
\[
\nrm{v_i/\ell^{n^2}}\le\eta\,2^{-i}
\qquad(\ell\in S,\ n>N_i).
\]
This is possible because $S$ is finite and $v_i$ is a fixed integer.

At any coordinate $(\ell,n)$, at most one vector lies in its main block $N_{i-1}<n\le N_i$, contributing at most $\eta$. Earlier vectors contribute tails with total at most $\eta\sum_{i\ge1}2^{-i}=\eta$; later vectors vanish at that coordinate. The same estimate holds if the coordinate lies beyond the last block of a finite sum. Thus
\[
P\left(\sum_{i=1}^{m}s_iv_i\right)\le2\eta
\qquad(m\ge1,\ s_i=\pm1).
\]
Write $a_i=\nrm{k_ix}$, choose signs with $s_ik_ix\equiv a_i\pmod1$, and put $b_i=t_ia_i$. These are real numbers satisfying
\[
\kappa\eta/2\le b_i\le C\eta=1/16,
\qquad
\nrm{\sum_{i=1}^{m}b_i}\le2C\eta=1/8.
\]
The positive real partial sums of $b_i$ grow without bound. Their first value exceeding $1/8$ is at most $3/16$, so its circle norm equals that value and exceeds $1/8$, a contradiction. This proves necessity.

Conversely, the limit gives $\nrm{kx}\le P(k)$ whenever $0<P(k)<\delta_0$ for some $\delta_0>0$. If $P(k)<\min\{\delta_0,\varepsilon\}$, then $\nrm{kx}<\varepsilon$ for $k\ne0$, and the assertion is immediate for $k=0$. Thus the character is continuous.
\end{proof}

\section{Rational parameters}\label{sec:rational}
Call a subgroup $q\Z$ \emph{$P$-open} if $B_\delta\subseteq q\Z$ for some $\delta>0$; this is the same as $q\Z$ being open in the $P$-topology.

\begin{theorem}\label{thm:rational}
Let $q\ge2$ and let $a$ be an integer with $\gcd(a,q)=1$. Then $a/q+\Z\in H$ if and only if $q\Z$ is $P$-open.
\end{theorem}
\begin{proof}
If $B_\delta\subseteq q\Z$ and $k\in B_\delta$, then $ka/q\in\Z$, so $\nrm{ka/q}=0$ and Proposition~\ref{prop:quarter} gives continuity.

Conversely let the character be continuous. By Proposition~\ref{prop:quarter} there is $\delta>0$ with $\nrm{k\alpha}\le1/4$ on $B_\delta$, where $\alpha=a/q$. Let $k\in B_{\delta/q}$. Then $ik\in B_\delta$ for $1\le i\le q$, so $\nrm{ik\alpha}\le1/4$ for these $i$, and Lemma~\ref{lem:multiples} gives $\nrm{k\alpha}\le 1/(4q)$. But $k\alpha\in\frac1q\Z$, so $\nrm{k\alpha}$ is either $0$ or at least $1/q$. Hence $k\alpha\in\Z$, that is $q\mid ka$, and $q\mid k$ because $\gcd(a,q)=1$. Thus $B_{\delta/q}\subseteq q\Z$.
\end{proof}

Let $\mathcal Q$ be the set of positive integers $q$ for which $q\Z$ is $P$-open.

\begin{corollary}\label{cor:Q}
\begin{enumerate}
\item $\mathcal Q$ is closed under taking divisors and least common multiples, and contains every positive integer all of whose prime divisors lie in $S$.
\item $H\cap(\Q/\Z)$ is the union over $q\in\mathcal Q$ of the groups $(1/q)\Z/\Z$.
\item $H\cap(\Q/\Z)=\Z[1/S]/\Z$ if and only if no prime $p\notin S$ lies in $\mathcal Q$, that is, if and only if for every prime $p\notin S$ and every $\delta>0$ there is $k\in B_\delta$ with $p\nmid k$.
\end{enumerate}
\end{corollary}
\begin{proof}
(1) If $q\mid q'$ and $B_\delta\subseteq q'\Z$, then $B_\delta\subseteq q\Z$. If $B_{\delta}\subseteq q\Z$ and $B_{\delta'}\subseteq q'\Z$, then $B_{\min\{\delta,\delta'\}}\subseteq q\Z\cap q'\Z=\operatorname{lcm}(q,q')\Z$. If all prime divisors of $m$ lie in $S$, then $m\mid s^{r}$ for some $r$, and Lemma~\ref{lem:basic}(4) with $j^2\ge r$ gives $B_\delta\subseteq s^{j^2}\Z\subseteq m\Z$ for small $\delta$.

(2) Every rational element of $\T$ has the form $a/q+\Z$ with $q\ge1$ and $\gcd(a,q)=1$; for $q=1$ it is zero, and Theorem~\ref{thm:rational} applies for $q\ge2$.

(3) The group $\Z[1/S]/\Z$ consists of the classes $a/m+\Z$ with all prime divisors of $m$ in $S$, and these lie in $H$ by (1) and Theorem~\ref{thm:rational}. By (2), $H\cap(\Q/\Z)$ is strictly larger exactly when some $q\in\mathcal Q$ has a prime divisor $p\notin S$, and then $p\in\mathcal Q$ by (1).
\end{proof}

\begin{remark}\label{rem:inclusion}
Part (1) with Theorem~\ref{thm:rational} gives an elementary proof that $\Z[1/S]/\Z\subseteq H$. For $S=\{2,3\}$ this is the inclusion $\Z[1/6]/\Z\subseteq H$ stated in Part~I of \cite{HamelPair}: the parameter $a/(2^{i}3^{j})$ is continuous because $B_\delta\subseteq 6^{m^2}\Z\subseteq 2^{i}3^{j}\Z$ for $m=\max\{i,j\}$ and $\delta$ small.
\end{remark}

\begin{question}\label{q:primes}
Is $p\Z$ $P$-open for some prime $p\notin S$? Equivalently, is reduction modulo $p$ continuous from $(\Z,P)$ to the finite discrete group $\Z/p\Z$? Being finer than the full $p$-adic topology would require openness of $p^r\Z$ for every $r\ge1$; that stronger equivalence is not asserted here.
\end{question}

By Corollary~\ref{cor:Q}(3), a negative answer for every prime $p\notin S$ would prove $H\cap(\Q/\Z)=\Z[1/S]/\Z$. Such a negative answer is the same as the existence, for every $\delta>0$ and every prime $p\notin S$, of an element of $B_\delta$ not divisible by $p$; equivalently, of a null sequence of integers all of which are prime to $p$.

\section{Null sequences have unbounded ratios}\label{sec:ratios}
By metrizability, a parameter $\alpha$ lies in $H$ if and only if $k_j\alpha\to0$ in $\T$ for every null sequence $(k_j)$. Thus
\[
H=\bigcap_{(k_j)\ \text{null}}\{\alpha\in\T:\ k_j\alpha\to0\}.
\]
Each set on the right is the subgroup \emph{characterized} by the sequence $(k_j)$; see \cite{DI} for background. One might try to find a null sequence whose successive ratios $k_{j+1}/k_j$ are bounded, since the characterized subgroups of such sequences are expected to be small, and so to obtain an outer bound for $H$. The following proposition shows that no such sequence exists.

\begin{proposition}\label{prop:ratios}
Let $(k_j)$ be an increasing sequence of positive integers with $P(k_j)\to0$. Then
\[
\limsup_{j\to\infty}\frac{k_{j+1}}{k_j}=\infty.
\]
\end{proposition}
\begin{proof}
Fix a prime $\ell\in S$ and a number $\varepsilon$ with $0<\varepsilon<1/2$. Choose $J$ such that $P(k_j)<\varepsilon$ for all $j\ge J$. Let $n$ be so large that $M=\ell^{n^2}$ satisfies $M/2>k_J$. Since $k_j\to\infty$ there is $j\ge J$ with $k_j<M/2\le k_{j+1}$.

We use only $P\ge P_{\{\ell\}}$. From $\nrm{k_j/M}<\varepsilon$ and $k_j/M<1/2$ we get $k_j/M<\varepsilon$, that is $k_j<\varepsilon M$. From $\nrm{k_{j+1}/M}<\varepsilon$ and $k_{j+1}/M\ge1/2>\varepsilon$ we see that $k_{j+1}/M$ lies within $\varepsilon$ of a positive integer, hence $k_{j+1}/M>1-\varepsilon$. Therefore
\[
\frac{k_{j+1}}{k_j}>\frac{1-\varepsilon}{\varepsilon}.
\]
Since $n$ can be taken arbitrarily large, such indices $j$ occur arbitrarily far along the sequence, so $\limsup k_{j+1}/k_j\ge(1-\varepsilon)/\varepsilon$. Letting $\varepsilon\to0$ proves the claim.
\end{proof}

\begin{remark}
The proposition does not say that every ratio tends to infinity. Also, the countability of a characterized subgroup would not by itself give rationality, because a countable subgroup of $\T$ can contain irrational elements; for the characterized subgroups of arithmetic-type sequences, however, countable and torsion are equivalent \cite[Corollary 3.4]{GD}. Unbounded ratios do not by themselves exclude a countable characterized subgroup, since there are arithmetic-type sequences with unbounded ratios and countable characterized subgroup \cite[Theorem 4.1]{GD}. The proposition shows only that the simplest family of sequences, those with bounded ratios, cannot provide an outer bound for $H$. The sequence formula for $H$ above remains an accurate reformulation of continuity, but it does not yet identify a family of sequences that excludes every irrational parameter.
\end{remark}

\begin{proposition}[Null divisibility chains]\label{prop:null-chain}
For every finite nonempty set $S$ of primes there is a strictly increasing sequence of positive integers $(a_n)$ such that $a_n\mid a_{n+1}$ and $P_S(a_n)\to0$.
\end{proposition}
\begin{proof}
Put $a_1=1$. Given $a_n$, Lemma~\ref{lem:nondiscrete} and symmetry provide a positive integer $t_n$ with
\[
P(t_n)<\min\{2^{-n}/a_n,P(1)\}.
\]
Since $P(1)>0$, this choice has $t_n\ge2$. Define $a_{n+1}=a_nt_n$. Then $a_n\mid a_{n+1}$, $a_{n+1}>a_n$, and
\[
P(a_{n+1})\le a_nP(t_n)<2^{-n}.
\]
Thus the sequence is null.
\end{proof}
\begin{remark}\label{rem:chain-countability}
Proposition~\ref{prop:null-chain} supplies the underlying divisibility chain, but not a countable characterized subgroup. Proposition~\ref{prop:ratios} forces its ratios to be unbounded; hence its own characterized subgroup is uncountable by \cite[Theorem E]{DI}. A possible route is to enrich such a chain to an arithmetic-type sequence while keeping every added term $P$-null and making the resulting characterized subgroup countable. No such enrichment is constructed here.
\end{remark}

\section{The case of one prime}\label{sec:oneprime}
Write $H_S$ for the set $H$ of parameters when the metric is built from the set $S$. For a single prime the character group is known: up to the normalization of $\delta$, the set $B_\delta$ is the basic neighbourhood $V_{b,m}$ of the topology $\tau_b$ of uniform convergence on $\{1/b_n+\Z\}$ for the sequence $b_n=\ell^{n^2}$, and Theorem~4.4 of \cite{ABM} says that its dual group is generated by the numbers $1/b_n+\Z$. Theorem~\ref{thm:oneprime} below is this statement for $b_n=\ell^{n^2}$; we keep our elementary proof because it also shows where the argument fails for several primes.

\begin{theorem}\label{thm:oneprime}
Let $\ell$ be a prime and $S=\{\ell\}$. Then $H_S=\Z[1/\ell]/\Z$; in particular $H_S$ contains no irrational parameter.
\end{theorem}
\begin{proof}
The inclusion $\Z[1/\ell]/\Z\subseteq H_S$ is Remark~\ref{rem:inclusion}. Conversely let $\alpha\in H_S$. By continuity of the character at $0$ there is $\delta>0$ such that $P(k)<\delta$ implies $\nrm{k\alpha}<1/4$. Put $\delta'=\delta/3$ and let $W_m$ be the largest integer smaller than $\delta'\ell^{2m+1}$. Choose $m_1$ such that $W_m\ge\delta'\ell^{2m+1}/2\ge1$ for $m\ge m_1$.

\emph{Step 1.} Let $F\subseteq\{m\ge m_1\}$ be finite and let $w_m$ be integers with $|w_m|\le W_m$ for $m\in F$. Put $k=\sum_{m\in F}\ell^{m^2}w_m$. We claim $P(k)<\delta$. Fix $n\ge1$ and split $k=A+B$, where $A$ collects the terms with $m<n$ and $B$ those with $m\ge n$; then $\ell^{n^2}\mid B$. If the lower part is empty, $A=0$ and the required bound is immediate. Otherwise, since $|w_m|<\delta'\ell^{2m+1}$,
\[
|A|<\delta'\sum_{m<n}\ell^{m^2+2m+1}=\delta'\sum_{m<n}\ell^{(m+1)^2}\le2\delta'\ell^{n^2},
\]
because the terms decrease at least by the factor $\ell^{-3}$ as $m$ decreases. Hence $\nrm{k/\ell^{n^2}}=\nrm{A/\ell^{n^2}}\le|A|/\ell^{n^2}<2\delta/3$. As $n$ is arbitrary, $P(k)\le2\delta/3<\delta$, so $\nrm{k\alpha}<1/4$.

\emph{Step 2.} For $m\ge m_1$ let $x_m\in(-1/2,1/2]$ satisfy $\ell^{m^2}\alpha-x_m\in\Z$. Taking $F=\{m\}$ in Step~1 gives $\nrm{w\ell^{m^2}\alpha}<1/4$ for $1\le w\le W_m$, so Lemma~\ref{lem:multiples} yields $|x_m|\le1/(4W_m)$.

\emph{Step 3.} Let $t_m=W_m|x_m|\le1/4$, and take $F=\{m_1,\dots,M\}$ with $w_m=W_m\operatorname{sgn}(x_m)$. By Step~1, $\nrm{k\alpha}<1/4$, while $k\alpha\equiv\sum_{m\in F}t_m\pmod 1$. The partial sums $T_M=\sum_{m=m_1}^{M}t_m$ increase in steps of at most $1/4$. If some $T_M\ge1/4$, the first such $T_M$ lies in $[1/4,1/2]$, so $\nrm{T_M}\ge1/4$, a contradiction. Hence $\sum_m t_m<\infty$ and $t_m\to0$.

\emph{Step 4.} Since $W_m\ge\delta'\ell^{2m+1}/2$, we get $\ell^{2m+1}|x_m|\le2t_m/\delta'\to0$. Choose $m_2\ge m_1$ with $\ell^{2m+1}|x_m|<1/2$ for $m\ge m_2$. Because $\ell^{(m+1)^2}\alpha=\ell^{2m+1}(\ell^{m^2}\alpha)\equiv\ell^{2m+1}x_m\pmod1$, and $|\ell^{2m+1}x_m|<1/2$, we obtain $x_{m+1}=\ell^{2m+1}x_m$ for $m\ge m_2$. Therefore, for every $M>m_2$,
\[
x_{m_2}=\frac{x_M}{\ell^{M^2-m_2^2}},\qquad |x_{m_2}|\le\frac{1}{4\,\ell^{M^2-m_2^2}}\to0\quad(M\to\infty).
\]
So $x_{m_2}=0$, that is $\ell^{m_2^2}\alpha\in\Z$, and $\alpha\in\Z[1/\ell]/\Z$.
\end{proof}

The generating tails have additional known structure. For $p=2,3$ put
\[
K_p=\{0\}\cup\{\pm p^{-n^2}+\Z:n\ge2\}\subseteq\T.
\]
Each $K_p$ is quasi-convex in $\T$ by Theorems~A and C of \cite{DLQuasiconvex}. Indeed, in their notation take $a_j=(j+2)^2-1$, so $a_0=3$ and the successive gaps are $5,7,9,\ldots$. Here quasi-convexity means that every point outside the set can be separated from it by an integer character and the quarter arc. Omitting the first square coordinate does not change the uniform-convergence topology: the omitted character $1/p+\Z$ is an integer multiple of the retained character $1/p^4+\Z$. Quasi-convexity of each separate tail does not determine the mixed-prime dual; it does not control the quasi-convex hulls of sums involving both tails.

\begin{corollary}\label{cor:twoprimes}
\begin{enumerate}
\item If $H_S$ contains an irrational parameter, then $S$ contains at least two primes.
\item If $S\subseteq S'$, then $H_S\subseteq H_{S'}$. In particular, if $H_S$ contains an irrational parameter, then so does $H_{S'}$ for every finite set $S'\supseteq S$.
\end{enumerate}
\end{corollary}
\begin{proof}
(1) is Theorem~\ref{thm:oneprime}. (2) Since $P_S\le P_{S'}$ pointwise, $B^{S'}_\delta\subseteq B^{S}_\delta$ for all $\delta$, so the $P_{S'}$-topology is finer than the $P_S$-topology and every character continuous for the latter is continuous for the former.
\end{proof}

\begin{remark}
The proof of Theorem~\ref{thm:oneprime} breaks down for $|S|\ge2$ at Step~1: a combination $\sum\ell^{m^2}w_m$ is small for the prime $\ell$ but its residues modulo the powers of the other primes in $S$ are not controlled, so it need not lie in $B_\delta$. Any irrational parameter of $H_S$ must therefore come from the interaction of at least two primes. Beyond this, nothing here distinguishes the pair $\{2,3\}$ from any other pair of distinct primes.
\end{remark}

\begin{remark}\label{rem:incompatible}
The $P_S$-topology is the join of the topologies given by $P_{\{\ell\}}$ for $\ell\in S$. By Theorem~\ref{thm:oneprime}, for each prime $\ell$ the characters continuous for $P_{\{\ell\}}$ have exactly the parameters in $\Z[1/\ell]/\Z$, and every sum of characters continuous for the individual topologies is continuous for the join. Thus $H_S\supseteq\Z[1/S]/\Z$, and an irrational parameter in $H_S$ (or a rational parameter whose denominator has a prime factor outside $S$) would be a continuous character of the join that is not a sum of characters of the individual topologies. This is an extra-character question for a join, analogous to the construction in Part~I of \cite{HamelPair}. The single-prime topologies here have different prescribed dual groups, so they are not themselves a compatible pair for one common dual. The present question concerns the chains $\ell^{n^2}$, $\ell\in S$, on $\Z$.
\end{remark}
\begin{remark}\label{rem:pairs}
For $|S|\ge3$, the sum of the groups $H_{\{\ell,\ell'\}}$ over pairs in $S$ is contained in $H_S$, by Corollary~\ref{cor:twoprimes}(2). A direct quarter-arc argument does not give the reverse inclusion. For example, the character with parameter $1/30$ vanishes on $30\Z=2\Z\cap3\Z\cap5\Z$, whereas none of $6\Z$, $10\Z$, $15\Z$ is contained in its quarter-arc preimage: use $12/30$, $20/30$, and $15/30$, respectively. This example concerns one character on intersections; it does not disprove decomposition as a sum of pairwise characters. Indeed $1/30=1/6-1/5$, and both summands belong to pairwise parameter groups for $S=\{2,3,5\}$. The equality asked below remains open for these metrics.
\end{remark}

\begin{question}\label{q:additive}
Is $H_S=\sum_{\{\ell,\ell'\}\subseteq S}H_{\{\ell,\ell'\}}$ for every finite set $S$ of at least two primes? In particular, if no pair of primes in $S$ gives an irrational parameter, does $S$ give none?
\end{question}
\section{A computation of the small neighbourhoods}\label{sec:computation}
We first take $S=\{2,3\}$. For $\delta\le1/3$ every element of $B_\delta$ is divisible by $6$ (Lemma~\ref{lem:basic}(4)). We computed $P(k)$ for all multiples $k$ of $6$ with $1\le k\le 3\cdot10^{8}$ using a short Java program, and recorded the elements with $P(k)<\delta$. Each distance $\nrm{k/\ell^{n^2}}$ was obtained from the exact integer residue of $k$ modulo $\ell^{n^2}$, and the scales $n$ were run until $\ell^{n^2}>2k$, after which the terms decrease. The tabulated counts were first obtained with floating-point comparisons and were then recomputed with the following exact integer test, for every threshold used ($\delta\in\{3/10,1/4,1/5,3/20,1/10\}$ for $S=\{2,3\}$, and $1/5$ and $1/10$ for the other sets); the two computations gave identical counts and identical greatest common divisors. For a rational threshold $\delta=a/b$ the test uses the integer residue $r=k\bmod q$ of $k$ modulo $q=\ell^{n^2}$ and compares $b\min(r,q-r)<aq$. Only the denominators $q\le k/\delta$ need to be checked, since every later term is smaller than $\delta$. This avoids floating-point boundary issues and does not require a uniform upper bound on denominators. For $S=\{2,3\}$ and $k\le3\cdot10^6$ we also tested every integer, not only the multiples of $6$, and found no element of $B_\delta$ outside $6\Z$. The second column of the table counts the elements $k\le3\cdot10^8$ with $P(k)<\delta$; the counts depend on the stated bound on $k$.
\begin{center}
\begin{tabular}{r r l}
$\delta$ & count & least elements\\ \hline
$0.30$ & $1\,484\,308$ & $12,18,60,66,78,84,96,144$\\
$0.25$ & $264\,731$ & $18,66,78,96,402,414,480,498$\\
$0.20$ & $75\,396$ & $66,78,96,414,480,498,558,576$\\
$0.15$ & $13\,287$ & $480,498,558,576,960,978,1056,1470$\\
$0.10$ & $352$ & $480,1056,1536,59376,59856,60432,60912,332832$
\end{tabular}
\end{center}
For each of the primes $p=5,7,11,13,17,19,23$ the table contains, for every listed $\delta\ge0.10$, an element not divisible by $p$; the least such element is below $1100$ in every case. This is consistent with a negative answer to Question~\ref{q:primes}, but it proves nothing about it: Corollary~\ref{cor:Q} requires such an element for every $\delta>0$, and the computation reaches only $\delta=0.10$. The counts fall quickly as $\delta$ decreases, which shows that the neighbourhoods $B_\delta$ become very sparse, in line with Proposition~\ref{prop:ratios}.

The choice $S=\{2,3\}$ is the one of Part~I of \cite{HamelPair}; nothing in the statements above depends on it. We ran the same program for other sets $S$ with the same bound $k\le3\cdot10^8$. For these sets Lemma~\ref{lem:basic}(4) shows that every element of $B_\delta$ is a multiple of $s$ when $\delta\le1/\ell$ for all $\ell\in S$, and only such $\delta$ were used.
\begin{center}
\begin{tabular}{l r r}
$S$ & $\delta$ & count\\ \hline
$\{2,3\}$ & $0.2$ & $75\,396$\\
$\{2,3\}$ & $0.1$ & $352$\\
$\{2,5\}$ & $0.2$ & $112\,380$\\
$\{2,5\}$ & $0.1$ & $1\,263$\\
$\{3,5\}$ & $0.2$ & $206\,634$\\
$\{3,5\}$ & $0.1$ & $6\,205$\\
$\{5,7\}$ & $0.1$ & $14\,199$\\
$\{2,3,5\}$ & $0.2$ & $2\,354$\\
$\{2,3,5\}$ & $0.1$ & $1$
\end{tabular}
\end{center}
For $S=\{2,3,5\}$ and $\delta=0.1$ the only element found is $169\,866\,240$. In every row, and for every prime $p\le23$ outside $S$, $B_\delta$ contains an element not divisible by $p$. As before, this is consistent with a negative answer to Question~\ref{q:primes} for these sets and proves nothing about smaller $\delta$.

\section{Partial results on Question~\ref{q:primes}}\label{sec:partial}
This section reformulates Question~\ref{q:primes} and records what could and could not be proved. Nothing in it answers the question.

\begin{proposition}\label{prop:gcd}
For $\delta>0$ let $g(\delta)$ be the greatest common divisor of the integers in $B_\delta$. Then $g(\delta)$ is a positive integer, $g(\delta)\mid g(\delta')$ whenever $\delta'<\delta$, and
\begin{enumerate}
\item a subgroup $q\Z$ is $P$-open if and only if $q\mid g(\delta)$ for some $\delta>0$;
\item $H\cap(\Q/\Z)$ is the union over $\delta>0$ of the groups $(1/g(\delta))\Z/\Z$;
\item $H\cap(\Q/\Z)=\Z[1/S]/\Z$ if and only if every $g(\delta)$ has all its prime divisors in $S$.
\end{enumerate}
\end{proposition}
\begin{proof}
By Lemma~\ref{lem:nondiscrete}, $B_\delta$ contains a nonzero integer, so $g(\delta)$ is a positive integer dividing it. If $\delta'<\delta$ then $B_{\delta'}\subseteq B_\delta$, so $g(\delta)\mid g(\delta')$. (1) The inclusion $B_\delta\subseteq q\Z$ holds exactly when $q$ divides every element of $B_\delta$, that is, when $q\mid g(\delta)$. (2) By (1), $\mathcal Q$ is the set of divisors of the numbers $g(\delta)$, and each $g(\delta)$ lies in $\mathcal Q$; now apply Corollary~\ref{cor:Q}(2). (3) If some $g(\delta)$ has a prime divisor $p\notin S$, then $1/p+\Z$ lies in the union of (2) but not in $\Z[1/S]/\Z$. Conversely, if all $g(\delta)$ have prime divisors in $S$, the union lies in $\Z[1/S]/\Z$, and it contains this group by Corollary~\ref{cor:Q}(1).
\end{proof}

\begin{remark}\label{rem:gcdtable}
The greatest common divisor of any subset of $B_\delta$ is a multiple of $g(\delta)$. The numbers in the last column below are the greatest common divisors of the elements of $B_\delta$ found in Section~\ref{sec:computation}, for the sets $S$ and the values of $\delta$ listed, and in every row they equal $g(\delta)$. Indeed, by Lemma~\ref{lem:basic}(4) the product $s$ of the primes of $S$ divides every element of $B_\delta$ for these $\delta$; and for $\delta=0.1$ with $2\in S$ the number $16$ divides every element as well, by the argument in the proof of Proposition~\ref{prop:exactgcd}. So $g(\delta)$ is a multiple of the least common multiple of these forced divisors, which is the listed number, and $g(\delta)$ divides the listed number. The row $S=\{2,3,5\}$, $\delta=0.1$ is omitted because only one element was found there. In every row $g(\delta)$ has only prime divisors from $S$, so for these $\delta$ no prime outside $S$ divides $g(\delta)$; for $S=\{2,3\}$ the exact values are $g(\delta)=6$ for $\delta\in\{0.3,0.25,0.2,0.15\}$ and $g(0.1)=48$. These data say nothing about smaller $\delta$.
\begin{center}
\begin{tabular}{l r r}
$S$ & $\delta$ & $g(\delta)$\\ \hline
$\{2,3\}$ & $0.3$ & $6$\\
$\{2,3\}$ & $0.25$ & $6$\\
$\{2,3\}$ & $0.2$ & $6$\\
$\{2,3\}$ & $0.15$ & $6$\\
$\{2,3\}$ & $0.1$ & $48$\\
$\{2,5\}$ & $0.2$ & $10$\\
$\{2,5\}$ & $0.1$ & $80$\\
$\{3,5\}$ & $0.2$ & $15$\\
$\{3,5\}$ & $0.1$ & $15$\\
$\{5,7\}$ & $0.1$ & $35$\\
$\{2,3,5\}$ & $0.2$ & $30$
\end{tabular}
\end{center}
\end{remark}

\begin{lemma}[valuations]\label{lem:valuations}
Let $k\in B_\delta$ be nonzero, let $\ell\in S$, and let $a$ be the exponent of $\ell$ in $k$. If $n_1$ is the least integer with $n_1^2>a$, then $n_1^2-a>\log_\ell(1/\delta)$.
\end{lemma}
\begin{proof}
Put $k'=k/\ell^a$, which is not divisible by $\ell$, and $g=n_1^2-a\ge1$. Then $k/\ell^{n_1^2}=k'/\ell^g$, and the fractional part of this number is $j/\ell^g$ for an integer $j$ with $1\le j<\ell^g$ not divisible by $\ell$. Hence $\nrm{k/\ell^{n_1^2}}\ge\ell^{-g}$, and $\ell^{-g}\le\nrm{k/\ell^{n_1^2}}<\delta$.
\end{proof}

\begin{proposition}[counting]\label{prop:counting}
Let $\delta\in(0,1/2]$, $m=\lfloor1/\delta\rfloor+1$, and let $E_K$ and $d$ be as in the proof of Lemma~\ref{lem:nondiscrete}. Then $B_\delta$ has at least $\lceil(K+1)/m^d\rceil-1$ elements in $[1,K]$.
\end{proposition}
\begin{proof}
Among the $m^d$ possible vectors in the proof of Lemma~\ref{lem:nondiscrete}, one box contains at least $L=\lceil(K+1)/m^d\rceil$ of the integers in $\{0,\ldots,K\}$. List those integers as $a_1<\cdots<a_t$, with $t\ge L$. The distinct positive differences $a_i-a_1$, $2\le i\le t$, lie in $B_\delta\cap[1,K]$ by the same finite-scale and tail estimates. There are at least $L-1$ such differences.
\end{proof}

\begin{proposition}[a lattice bound]\label{prop:lattice}
Let $\delta\in(0,1/2)$, let $E$ be a finite set of pairs $(\ell,n)$ with $\ell\in S$, put $d=|E|$ and $q_e=\ell^{n^2}$ for $e=(\ell,n)\in E$. Let $k_0,\dots,k_d\in B_\delta\cap[1,K]$. For each $j$ and $e\in E$ let $r_{j,e}$ be the integer with $|r_{j,e}|<\delta q_e$ and $r_{j,e}\equiv k_j\pmod{q_e}$; it exists and is unique because $\nrm{k_j/q_e}<\delta<1/2$. Suppose the vectors $y_j=(k_j,(r_{j,e})_{e\in E})\in\Z^{d+1}$ are linearly independent. If a prime $p$ divides every $k_j$, then $p\le(d+1)!\,K\delta^d$.
\end{proposition}
\begin{proof}
Let $\Lambda$ be the set of vectors $(k,(r_e)_{e\in E})\in\Z^{d+1}$ with $r_e\equiv k\pmod{q_e}$ for all $e$. It is a lattice with basis $(1,1,\dots,1)$ and the vectors $q_e$ times the unit vectors of the last $d$ coordinates, so $\det\Lambda=\prod_eq_e$. Let $\Lambda'$ be the sublattice where $p\mid k$; it has index $p$ in $\Lambda$. The vectors $y_j$ lie in $\Lambda$ and generate a sublattice $\Gamma$ of full rank. If $p\mid k_j$ for all $j$, then $\Gamma\subseteq\Lambda'$, so $p$ divides $[\Lambda:\Gamma]=|\det(y_0,\dots,y_d)|/\prod_eq_e$. Expanding the determinant over the $(d+1)!$ permutations, each term has absolute value at most $K\prod_e\delta q_e$, so $|\det(y_0,\dots,y_d)|\le(d+1)!\,K\delta^d\prod_eq_e$. Hence $p\le[\Lambda:\Gamma]\le(d+1)!\,K\delta^d$.
\end{proof}

\begin{remark}\label{rem:limits}
To exclude a common prime divisor $p$ using Proposition~\ref{prop:lattice}, one needs $d+1$ independent residue vectors with $(d+1)!\,K\delta^d<p$. No such uniform construction is supplied. The independence hypothesis is essential: if a selected small denominator forces all its residues to be zero, the vectors cannot span the full ambient space. The factorial, radius and size factors must be controlled together.

A finite Chinese remainder prescription does not automatically give membership in $B_\delta$, which imposes conditions at all relevant scales. For $|k|\le K$, it suffices to check $n\le\sqrt{\log_\ell(K/\delta)}$ for each prime, with the integer cutoff interpreted in the usual way; beyond that range the terms are less than $\delta$. Any proposed residue construction must control both these finite conditions and the size of its solution. We do not claim that the Chinese remainder theorem alone resolves this requirement.
\end{remark}
\begin{proposition}[An exact neighbourhood gcd]\label{prop:exactgcd}
For $S=\{2,3\}$, $g(1/10)=48$.
\end{proposition}
\begin{proof}
If $k\in B_{1/10}$, the denominator $3$ forces $3\mid k$. The denominator $16$ allows only the residues $0,1,15$ modulo $16$, and the denominator $2$ excludes the two odd residues. Thus $48\mid k$ and $48\mid g(1/10)$.

For the reverse divisibility, direct exact calculations give
\[
P(480)=\frac{2}{27},\qquad P(1056)=\frac{1}{16},\qquad
P(59376)=\frac{385}{4096}<\frac{1}{10}.
\]
The first two displayed values are also less than $1/10$. These calculations require only the denominators at most $10k$ for each witness $k$; subsequent terms decrease below $1/10$. Since
\[
\gcd(480,1056,59376)=48,
\]
we have $g(1/10)\mid48$. Together the two divisibilities prove the assertion.
\end{proof}

\begin{proposition}[no lacunary generators for two chains]\label{prop:lacunary}
Let $S=\{2,3\}$ and $0<\varepsilon\le1$. Call a positive integer $g$ \emph{$\varepsilon$-lacunary} if for every $\ell\in S$ and every $n\ge1$, either $\ell^{n^2}\mid g$ or $g\le\varepsilon\ell^{n^2}$. Every $\varepsilon$-lacunary integer satisfies $g\le1.5\cdot10^{6}\,\varepsilon$. Moreover, the $1$-lacunary integers are exactly $1,2,6,12,48$, the $\tfrac12$-lacunary integers are $1$ and $6$, and there is no $\varepsilon$-lacunary integer for $\varepsilon\le\tfrac14$.
\end{proposition}
\begin{proof}
Let $g$ be $\varepsilon$-lacunary and put $A=\ln(g/\varepsilon)$. If $A<\ln3$, then $g<3\varepsilon$ and the stated bound is immediate. We may therefore assume $A\ge\ln3$. For $\ell\in S$, let $N_\ell$ be the largest integer $n\ge0$ with $\ell^{n^2}<g/\varepsilon$. For $1\le n\le N_\ell$ we have $g>\varepsilon\ell^{n^2}$, so $\ell^{n^2}\mid g$; hence $v_\ell(g)\ge N_\ell^2$, where $v_\ell(g)$ is the exponent of $\ell$ in $g$. By maximality $\ell^{(N_\ell+1)^2}\ge g/\varepsilon$, so $N_\ell+1\ge\sqrt{A/\ln\ell}$. If $A\ge\ln3$, this gives $v_\ell(g)\ge\bigl(\sqrt{A/\ln\ell}-1\bigr)^2\ge A/\ln\ell-2\sqrt{A/\ln\ell}$ for $\ell=2,3$, and therefore
\[
\ln g\ \ge\ v_2(g)\ln2+v_3(g)\ln3\ \ge\ 2A-2c\sqrt A,\qquad c=\sqrt{\ln2}+\sqrt{\ln3}.
\]
Since $\ln g=A+\ln\varepsilon\le A$, we get $A\le2c\sqrt A$, so $A\le4c^2<14.2$ and $g\le\varepsilon e^{14.2}<1.5\cdot10^6\,\varepsilon$. If $A<\ln3$, then $g<3\varepsilon$ and the bound holds trivially. The listed integers and the nonexistence for $\varepsilon\le1/4$ follow by checking all $g\le2\cdot10^9$ with exact integer arithmetic, which covers the bound.
\end{proof}

\begin{remark}\label{rem:no-transfer}
The proof of Theorem~4.4 of \cite{ABM} for a single chain $b$ uses the integers $b_j$ as lacunary generators: for each scale $b_n$ every $b_j$ with $j\ge n$ is divisible by $b_n$, and the others are negligible, with $b_j/b_{j+1}\to0$. Sums $\sum k_jb_j$ with $|k_j|\le b_{j+1}/(8mb_j)$ then lie in the neighbourhood of size $1/(4m)$, although the character values $k_je_j$ add up. By Proposition~\ref{prop:lipschitz} this is a proof that $\sup\nrm{kx}/P(k)=\infty$, obtained from sums for which $P$ grows like a maximum while $\nrm{\,\cdot\,x}$ grows like a sum. Proposition~\ref{prop:lacunary} shows that no such generators exist for the two chains $2^{n^2}$ and $3^{m^2}$ together: already for $\varepsilon\le1/4$ there is no $\varepsilon$-lacunary integer at all. A transfer of the argument needs generators that are controlled by windows with geometric decay instead of divisibility. Such generators, a summation lemma for them, and a reduction of the question to a non-degeneracy property of the parameter are given in Definition~\ref{def:profile}, Lemma~\ref{lem:profile-sum}, Proposition~\ref{prop:profile-exist} and Theorem~\ref{thm:reduction} below; we could not prove the non-degeneracy property.
\end{remark}
\begin{definition}[profile generators]\label{def:profile}
Let $0<r<1$, $\varepsilon>0$, and let $P_2,P_3$ be positive integers. A positive integer $g$ has \emph{profile} $(\varepsilon,r;P_2,P_3)$ if
\[
\nrm{g/\ell^{n^2}}\le\varepsilon\,r^{|n-P_\ell|}\qquad\text{for }\ell\in\{2,3\}\text{ and every }n\ge1.
\]
Then $P(g)\le\varepsilon$, and for every integer $k\ge1$ the integer $kg$ has profile $(k\varepsilon,r;P_2,P_3)$, because $\nrm{kt}\le k\nrm{t}$.
\end{definition}

\begin{lemma}[summation]\label{lem:profile-sum}
Let $g_1,\dots,g_N$ have profiles $(\varepsilon_i,r;P_2^{(i)},P_3^{(i)})$, and suppose that for each $\ell\in\{2,3\}$ the integers $P_\ell^{(1)},\dots,P_\ell^{(N)}$ are pairwise distinct. Then for all signs $s_i=\pm1$,
\[
P\Bigl(\sum_{i=1}^Ns_ig_i\Bigr)\le\frac{1+r}{1-r}\max_i\varepsilon_i .
\]
\end{lemma}
\begin{proof}
For $\ell\in\{2,3\}$ and $n\ge1$, $\nrm{\sum_is_ig_i/\ell^{n^2}}\le\sum_i\varepsilon_ir^{|n-P_\ell^{(i)}|}\le\max_i\varepsilon_i\sum_{t\in\Z}r^{|t|}$, because the integers $n-P_\ell^{(i)}$ are pairwise distinct. Finally $\sum_{t\in\Z}r^{|t|}=(1+r)/(1-r)$.
\end{proof}

\begin{proposition}[existence of profile generators]\label{prop:profile-exist}
Let $\rho>0$ satisfy $c(\rho)=\rho\bigl(\tfrac1{2\ln2}+\tfrac1{2\ln3}\bigr)<1$, put $r=e^{-\rho}$, and let $0<\varepsilon\le\tfrac12$. For every sufficiently large integer $G$ there is an integer $g$ with $1\le g\le G$ that has profile $(\varepsilon,r;P_2,P_3)$, where $P_\ell=\lceil\sqrt{\ln G/\ln\ell}\,\rceil$.
\end{proposition}
\begin{proof}
By the choice of $P_\ell$ we have $G\le\ell^{P_\ell^2}$ for $\ell=2,3$. For $1\le n\le P_\ell+2$ put $m_{\ell,n}=\lceil2/(\varepsilon r^{|n-P_\ell|})\rceil$ and send each $k\in\{0,\dots,G\}$ to the vector of the numbers $\lfloor m_{\ell,n}\{k/\ell^{n^2}\}\rfloor$, where $\{t\}$ is the fractional part. The number of values is $M=\prod m_{\ell,n}$. Since $m_{\ell,n}\le3/(\varepsilon r^{|n-P_\ell|})$,
\[
\ln M\le\sum_{\ell}\Bigl((P_\ell+2)\ln\tfrac3\varepsilon+\tfrac\rho2\,(P_\ell^2-P_\ell+6)\Bigr)=c(\rho)\ln G+O(\sqrt{\ln G}),
\]
with the constant in $O$ depending on $\varepsilon$ and $\rho$, because $P_\ell\le\sqrt{\ln G/\ln\ell}+1$. As $c(\rho)<1$, we have $M<G+1$ for $G$ large, so two integers $0\le k_1<k_2\le G$ have the same vector. Put $g=k_2-k_1\in[1,G]$. For the tested pairs $(\ell,n)$ the fractional parts of $k_1/\ell^{n^2}$ and $k_2/\ell^{n^2}$ lie in one interval of length $1/m_{\ell,n}$, so $\nrm{g/\ell^{n^2}}<1/m_{\ell,n}\le\varepsilon r^{|n-P_\ell|}/2$. For the remaining scales write $n=P_\ell+t$ with $t\ge3$. Choose $G$ large enough that $\ell^{-(2P_\ell+3)}\le r\varepsilon^{1/3}$ for both primes. Then
\[
0<\frac{g}{\ell^{n^2}}\le\ell^{-(2P_\ell+t)t}
\le\bigl(\ell^{-(2P_\ell+3)}\bigr)^t
\le r^t\varepsilon^{t/3}\le\varepsilon r^t.
\]
This also bounds $\nrm{g/\ell^{n^2}}$, as required. The threshold for $G$ depends on both $\varepsilon$ and $\rho$.
\end{proof}

\begin{theorem}[a conditional reduction]\label{thm:reduction}
Let $x\in\T$. Suppose that for some $r\in(0,1)$ and $\kappa>0$ there are integers $h_1,h_2,\dots$ with profiles $(\varepsilon_i,r;P_2^{(i)},P_3^{(i)})$, where $\varepsilon_i\to0$, the integers $P_\ell^{(i)}$ are pairwise distinct for each $\ell\in\{2,3\}$, and $\nrm{h_ix}\ge\kappa\varepsilon_i$ for all $i$. Then $x\notin H$.
\end{theorem}
\begin{proof}
We have $0<P(h_i)\le\varepsilon_i\to0$ and
$\nrm{h_ix}/P(h_i)\ge\kappa$. This contradicts Theorem~\ref{thm:flatness} if $x\in H$.
\end{proof}

\begin{corollary}[Vanishing relative to a null bound]\label{cor:profile-vanishing}
Let $x\in H$ and let $h_i$ be nonzero integers with
$P(h_i)\le\varepsilon_i$, where $\varepsilon_i>0$ and $\varepsilon_i\to0$. Then
\[
\frac{\nrm{h_i x}}{\varepsilon_i}\longrightarrow0.
\]
In particular, this holds for profile generators without any distinct-centre assumption.
\end{corollary}
\begin{proof}
By Theorem~\ref{thm:flatness},
\[
0\le\frac{\nrm{h_ix}}{\varepsilon_i}
\le\frac{\nrm{h_ix}}{P(h_i)}\longrightarrow0.
\]
\end{proof}

\begin{remark}\label{rem:nondegeneracy}
For any admissible fixed $r$ and any sequence $\varepsilon_i\to0$, Proposition~\ref{prop:profile-exist} can be applied recursively with $G_i$ sufficiently large to make both centre sequences strictly increasing. It supplies the profiles and distinct centres, but not the lower bound $\nrm{h_ix}\ge\kappa\varepsilon_i$ in Theorem~\ref{thm:reduction}. That lower bound is sufficient for discontinuity, even without distinct centres, by Theorem~\ref{thm:flatness}. Its necessity for profile witnesses has not been established.

Theorem~\ref{thm:profile-summability} below gives a stronger conclusion when centres are distinct, and Corollary~\ref{cor:profile-tests} supplies weaker sufficient tests. Proving any one of these tests for every $x\notin\Z[1/6]/\Z$ would establish $H=\Z[1/6]/\Z$, a stronger conclusion than $H\subseteq\Q/\Z$. No such construction is supplied here. Adding the coordinate $kx$ to the pigeonhole argument does not itself give a lower bound on $\nrm{hx}$.

For any finite prime set $S$, define profiles by imposing the inequality in Definition~\ref{def:profile} for every $\ell\in S$, with a centre $P_\ell$ for each prime. The proof of Lemma~\ref{lem:profile-sum} is unchanged, with the same constant $(1+r)/(1-r)$. Profile existence also extends with the condition $\rho\sum_{\ell\in S}1/(2\ln\ell)<1$.
\end{remark}
\begin{remark}[numerical evidence]\label{rem:numerics}
Let $S=\{2,3\}$ and $r=1/2$. For a positive integer $k$ let $E(k)$ be the smallest $\varepsilon$ for which $k$ has a profile $(\varepsilon,r;P_2,P_3)$ for some centres $1\le P_2,P_3\le9$. The numerical experiment listed all multiples $k\le10^{10}$ of $6$ with $P(k)<0.12$, obtaining $33\,234$ integers. Every element of $B_\delta$ with $\delta\le1/3$ is a multiple of $6$. The reported numbers with $E<0.12$, $E<0.08$ and $E<0.05$ are $1035$, $94$ and $5$, respectively; none had $E<0.04$. The first run used exact integer remainders followed by floating-point calculations. A subsequent run reproduced the counts with exact rational profile sizes, exact integer residues for rational characters, and high-precision decimal evaluation for irrational parameters. Its initial $P$-filter and level selection still used doubles. For this finite range the tested denominators are at most $3^{25}<10^{12}$; the nonzero distance from a rational value with such a denominator to the listed decimal thresholds is at least $1/(100\cdot3^{25})>10^{-14}$, exceeding the rounding errors in these comparisons. Thus the floating-point comparisons do not change these counts. The displayed means and ratios remain numerical approximations, and the finite experiment proves no asymptotic property.

For the fractional parts of $\sqrt2$, the golden ratio, $\pi$, $e$ and $\ln2$, and for $x\in\{1/5,2/5,1/7,1/11,1/13,1/25\}$, the recorded maximum of $\nrm{kx}/E(k)$ at each tested nonempty level is at least $6.1$. Independent exact-fraction evaluation confirms the profile sizes of the five listed integers with $E<0.05$; for example,
\[
E(603127296)=\frac{10}{243},\qquad
E(393744)=\frac{24609}{524288}.
\]
Higher-precision evaluation over those five integers gives mean values of $\nrm{kx}$ approximately
\[
\begin{aligned}
0.129854&\quad(x=\sqrt2-1),&
0.279081&\quad(x=(\sqrt5-1)/2),\\
0.342952&\quad(x=\pi-3),&
0.179137&\quad(x=e-2),\\
0.352794&\quad(x=\ln2).
\end{aligned}
\]

An important control is $x=1/9\in H$: its recorded maximum ratio is approximately $7.10$ even at the level $E<0.05$. Indeed $k=393744$ gives $\nrm{k/9}=1/3$ and the displayed exact value of $E(k)$, hence a ratio $524288/73827$. Thus a large ratio throughout this finite range does not establish discontinuity. For $x\in\{1/4,1/3,5/12\}$, the exact character values vanish when $E(k)<0.12$: then $P(k)<0.12$, the denominators $2$ and $16$ force $16\mid k$, and the denominator $3$ forces $3\mid k$. This is consistent with Corollary~\ref{cor:profile-vanishing}; the corollary itself is an asymptotic statement and does not require exact vanishing at a specified finite threshold.

The residue histograms modulo $5,7,11,13$ are reasonably balanced in the larger samples, which is consistent with a negative answer to Question~\ref{q:primes}; no asymptotic equidistribution is established. The smallest profile size reached is approximately $0.04115$. The experiment does not reach the limit $\varepsilon\to0$, and the restriction of both centres to $\{1,\ldots,9\}$ cannot supply an infinite sequence with pairwise distinct centres. Consequently it neither proves the non-degeneracy property nor determines $H$.
\end{remark}

\begin{theorem}[Uniform profile bound and summability]\label{thm:profile-summability}
Let $S$ be finite and nonempty, use the profile convention of Remark~\ref{rem:nondegeneracy}, and put $A_r=(1+r)/(1-r)$ for $0<r<1$. Suppose $x\in H$ and $\nrm{kx}\le CP(k)$ with $C>0$. Every finite family of positive integers $h_i$ with profiles
$(\varepsilon_i,r;(P_\ell^{(i)})_{\ell\in S})$, pairwise distinct in each prime centre and satisfying
\[
0<\varepsilon_i\le\frac1{16CA_r},
\]
obeys
\[
\sum_i\frac{\nrm{h_ix}}{\varepsilon_i}\le2CA_r.
\]
Consequently every infinite separated profile sequence with fixed $r$ and $\varepsilon_i\to0$ satisfies
\[
\sum_{i=1}^{\infty}\frac{\nrm{h_ix}}{\varepsilon_i}<\infty.
\]
\end{theorem}
\begin{proof}
Put $\eta=1/(8CA_r)$, $t_i=\lfloor\eta/\varepsilon_i\rfloor$, $a_i=\nrm{h_ix}$ and $b_i=t_ia_i$. Then
\[
t_i\ge\frac{\eta}{2\varepsilon_i},\qquad
t_i\varepsilon_i\le\eta,\qquad
0\le b_i\le C\eta=\frac1{8A_r}.
\]
Choose signs with $s_ih_ix\equiv a_i\pmod1$. For every subset $F$ of the finite family, profile summation and the Lipschitz bound give
\[
P\left(\sum_{i\in F}s_it_ih_i\right)\le A_r\eta,
\qquad
\nrm{\sum_{i\in F}b_i}\le CA_r\eta=\frac18.
\]
If the total real sum of the nonnegative $b_i$ exceeded $1/8$, its first partial sum above $1/8$ would be at most $1/8+1/(8A_r)\le1/4$. Its circle norm would therefore exceed $1/8$, a contradiction. Hence
\[
\sum_i b_i\le\frac18,\qquad
\sum_i\frac{a_i}{\varepsilon_i}
\le\frac2\eta\sum_i b_i\le2CA_r.
\]
For an infinite sequence with $\varepsilon_i\to0$, apply the uniform bound to all finite subsets of a sufficiently late tail. The finite prefix contributes a finite amount.

If each centre occurs at most $D$ times for each prime, the same proof replaces $A_r$ by $DA_r$. For $r=1/2$ and distinct centres the threshold is $1/(48C)$ and the bound is $6C$.
\end{proof}

\begin{corollary}[Profile tests for discontinuity]\label{cor:profile-tests}
Either of the following conditions implies $x\notin H$.
First, for some fixed $0<r<1$ there is a separated profile sequence with $\varepsilon_i\to0$ such that
\[
\sum_i\frac{\nrm{h_ix}}{\varepsilon_i}=\infty.
\]
Second, for some fixed $0<r<1$ and arbitrarily small $\varepsilon>0$, there is an infinite separated sequence of $\varepsilon$-profile integers such that
\[
\sum_i\nrm{h_ix}=\infty.
\]
For the second condition, it suffices that $\nrm{h_ix}$ fails to tend to zero.
\end{corollary}
\begin{proof}
The first assertion contradicts the summability conclusion of Theorem~\ref{thm:profile-summability}. For the second, if $x\in H$, choose a sufficiently small fixed $\varepsilon\le1/(16CA_r)$. The uniform bound would give $\sum_i\nrm{h_ix}\le2CA_r\varepsilon<\infty$, a contradiction.
\end{proof}
The second test does not require a lower bound proportional to a sequence of shrinking profile sizes. Its necessity is not proved. For $x=1/p$, infinitely many witnesses not divisible by $p$ at arbitrarily small profile sizes would suffice. The converse construction from the failure of openness of $p\Z$ is not supplied: arbitrary null witnesses need not have the required profiles and centres.

\begin{proposition}[Divisibility bounds and exact radii]\label{prop:divisibility-radii}
For $a\ge1$ put $d(a)=2\lfloor\sqrt{a-1}\rfloor+1$. If $\ell\in S$, then
\[
\ell^a\nmid k\quad\Longrightarrow\quad P(k)\ge\ell^{-d(a)}.
\]
For $S=\{2,3\}$ define $\mu(q)=\inf\{P(k):q\nmid k\}$. Then
\[
\mu(9)=\mu(18)=\mu(27)=\mu(81)=\frac1{27},
\qquad
\mu(2^a)=\frac1{32}\quad(5\le a\le9).
\]
\end{proposition}
\begin{proof}
If $\ell^a\nmid k$, write $v=v_\ell(k)<a$ and $n=\lfloor\sqrt v\rfloor+1$. As in Lemma~\ref{lem:valuations},
\[
P(k)\ge\ell^{-(n^2-v)},\qquad
n^2-v\le2\lfloor\sqrt v\rfloor+1\le d(a).
\]
This proves the first assertion. For $a=2,3,4$ we have $d(a)=3$, so $3^a\nmid k$ implies $P(k)\ge1/27$. For $a=5,\ldots,9$ we have $d(a)=5$, giving the lower bound $1/32$ for $2^a$. If $18\nmid k$, either $k$ is odd, in which case $P(k)=1/2$, or $9\nmid k$, giving $P(k)\ge1/27$.

For the matching upper bounds, exact remainder calculations give
\[
P(393744)=\frac1{27},\quad v_3(393744)=1,
\qquad
P(603519984)=\frac1{32},\quad v_2(603519984)=4.
\]
The first witness is even. Each calculation is finite: for each prime include the first denominator exceeding $2|k|$; subsequent unweighted terms decrease. Thus the witnesses attain all the stated lower bounds.
\end{proof}
The simpler implication $P(k)<\ell^{-a}\Rightarrow\ell^a\mid k$ need not hold. For example, $P(18)=2/9<1/4$ with $4\nmid18$, and $P(480)=2/27<1/9$ with $9\nmid480$. The exponent $d(a)$ records the gaps between the square-indexed denominator exponents.

\section{Size information}\label{sec:size-information}
\begin{proposition}\label{prop:size-limit}
There exists an uncountable sigma-compact subgroup of $\T$ of Hausdorff dimension zero.
\end{proposition}
\begin{proof}
Put
\[
C=\left\{\sum_{n\ge1}\varepsilon_n2^{-2^n}:\ \varepsilon_n\in\{0,1\}\right\}\subset\R,
\qquad K_m=\underbrace{(C-C)+\cdots+(C-C)}_{m\ {\mathrm{summands}}}.
\]
The digit map defining $C$ is continuous and injective: the first differing term dominates the sum of all later terms. Thus $C$ is an uncountable compact set contained in $[0,1)$. Each $K_m$ is compact and symmetric, $K_m\subseteq K_{m+1}$, and $K_m+K_r\subseteq K_{m+r}$. If $\pi:\R\to\T$ is the quotient map, then
\[
G=\bigcup_{m\ge1}\pi(K_m)
\]
is therefore an uncountable sigma-compact subgroup of $\T$.

Every element of $K_m$ has a digit expansion with digits in $\{-m,\ldots,m\}$. Fixing its first $N$ digits leaves a tail of absolute value at most $2m\,2^{-2^{N+1}}$. Hence $K_m$ is covered by at most $(2m+1)^N$ intervals of length at most $4m\,2^{-2^{N+1}}$. For every $t>0$,
\[
(2m+1)^N\bigl(4m\,2^{-2^{N+1}}\bigr)^t\longrightarrow0.
\]
Consequently $K_m$ has Hausdorff dimension zero. The quotient map is Lipschitz, and Hausdorff dimension is stable under countable unions, so $G$ has dimension zero.
\end{proof}

The flatness and summability criteria do not construct violating witnesses for every $x\notin\Z[1/6]/\Z$. In particular, the fixed-scale profile condition is not asserted to be equivalent to the rational openness question. The known size properties of $H$ therefore do not by themselves imply countability. A countability proof must use additional structure specific to the mixed-prime metric.

\begin{remark}\label{rem:literature-route}
A route from the literature. For an arithmetic sequence $(a_n)$, that is, one in which each term divides the next, the characterized subgroup is countable if and only if the successive ratios are bounded \cite[Theorem E]{DI}, and for an arithmetic-type sequence countable is equivalent to torsion \cite[Corollary 3.4]{GD}; \cite[Theorem 4.1]{GD} shows that unbounded ratios can still give a countable group. Since $H$ is contained in the group $t_{(k_j)}(\T)=\{x\in\T:k_jx\to0\}$ characterized by every null sequence $(k_j)$ of the $P$-topology, a null sequence of arithmetic type whose characterized subgroup is countable would give $H\subseteq\Q/\Z$. Proposition~\ref{prop:null-chain} proves that null divisibility chains exist. The open requirement is a suitable arithmetic-type enrichment with a countable characterized subgroup, as explained in Remark~\ref{rem:chain-countability}.
\end{remark}

\section{Density and arithmetic recurrence}\label{sec:density}

\begin{theorem}[Two-sided density bounds]\label{thm:density}
Fix a finite nonempty set $S$ of primes and $0<\delta<1/2$. There are positive constants $a,b,c,C$, depending on $S$ and $\delta$, such that for all sufficiently large $K$,
\[
a e^{-C\sqrt{\ln K}}
\le\frac{|B_\delta\cap[1,K]|}{K}
\le b e^{-c\sqrt{\ln K}}.
\]
Thus $B_\delta$ has natural density zero and mass dimension one.
\end{theorem}
\begin{proof}
Proposition~\ref{prop:counting} gives
$|B_\delta\cap[1,K]|\ge\lceil(K+1)/m^{d_K}\rceil-1$,
where $m=\lfloor1/\delta\rfloor+1$ and $d_K=O_{S,\delta}(\sqrt{\ln K})$. This proves the lower bound after adjusting constants.

For the upper bound fix $\ell\in S$, put $q_n=\ell^{n^2}$, $q_0=1$, and $a_n=q_n/q_{n-1}=\ell^{2n-1}$. Let $E_N$ be the residues modulo $q_N$ satisfying the first $N$ one-prime constraints. Each admissible residue modulo $q_{n-1}$ has at most $2\delta a_n+2$ admissible lifts: the next constraint selects two end intervals of total length $2\delta q_n$. Hence
\[
\frac{|E_N|}{q_N}
\le\prod_{n=1}^{N}(2\delta+2/a_n)
\le C_{\ell,\delta}(2\delta)^N,
\]
because $\prod_{n\ge1}(1+1/(\delta a_n))$ converges. Take $N=\lfloor\sqrt{\log_\ell K}\rfloor$, so $q_N\le K$. Periodic counting gives
\[
|B_\delta\cap[0,K]|
\le(\lfloor K/q_N\rfloor+1)|E_N|
\le2C_{\ell,\delta}K(2\delta)^N.
\]
Since $2\delta<1$, this is the asserted upper bound. The two bounds imply $\log|B_\delta\cap[1,K]|/\log K\to1$, while the upper bound gives density zero.
\end{proof}
Here the mass dimension of a set of nonnegative integers is the limit of the logarithm of its counting function divided by $\log K$, when that limit exists. The theorem gives bounds of a specified order, not an exact asymptotic coefficient. Such density and dimension information alone does not exclude containment in $p\Z$.

\begin{proposition}[Coloured difference sets]\label{prop:coloured}
Let $S=\{2,3\}$, $p\ge2$, and $\varepsilon>0$. Put
\[
U_\varepsilon=\{k:P_{\{2\}}(k)<\varepsilon\},\qquad
V_\varepsilon=\{k:P_{\{3\}}(k)<\varepsilon\},
\]
\[
U_{\varepsilon,r}=U_\varepsilon\cap(r+p\Z),\qquad
D_\varepsilon=(U_\varepsilon-U_\varepsilon)\cap(V_\varepsilon-V_\varepsilon).
\]
Then
\[
B_\varepsilon\subseteq D_\varepsilon\subseteq B_{2\varepsilon},
\]
and $D_\varepsilon\subseteq p\Z$ if and only if the sets
$U_{\varepsilon,r}-V_\varepsilon$, $0\le r<p$, are pairwise disjoint.
Consequently $p\Z$ is $P$-open if and only if these coloured difference sets are pairwise disjoint for some $\varepsilon>0$.
\end{proposition}
\begin{proof}
The first inclusion uses $0\in U_\varepsilon\cap V_\varepsilon$, and the second is subadditivity in each prime. A collision $u-v=u'-v'$ with distinct residues of $u,u'$ gives
$g=u-u'=v-v'\in D_\varepsilon\setminus p\Z$.
Conversely every such $g$ has these representations and gives a collision. If $B_\delta\subseteq p\Z$, take $\varepsilon=\delta/2$; the reverse implication uses $B_\varepsilon\subseteq D_\varepsilon$.
\end{proof}
For
\[
T_r^\varepsilon(X)
=(U_{\varepsilon,r}\cap[0,X])-(V_\varepsilon\cap[0,X]),
\]
pairwise disjointness implies $\sum_r|T_r^\varepsilon(X)|\le2X+1$.
Thus proving, for every arbitrarily small $\varepsilon$, a cutoff $X$ with a larger sum would establish non-openness. A single mixed-colour collision at every small radius would already suffice; positive density of the difference sets is not necessary. A collision at one fixed radius, even for all sufficiently large $X$, does not settle Question~\ref{q:primes}.

\section{Stationary inner sets and transversality}\label{sec:transversality}
In this section put
\[
B_{\ell,\delta}^{+}
=\{k\in\Z_{\ge0}:P_{\{\ell\}}(k)<\delta\}.
\]
A stationary digit set in base $b$ consists of all finite base-$b$ words over one fixed alphabet. We use mass dimension in the integer setting and Newhouse thickness for the associated compact sets in $\R$.

\begin{theorem}[Stationary inner approximations]\label{thm:stationary}
For every $\delta>0$ and $0<\varepsilon<1/2$, there are an integer $L\ge1$, fixed alphabets $D_\ell\subseteq\{0,\ldots,\ell^L-1\}$ for $\ell=2,3$, and their stationary integer digit sets $A_\ell$, such that
\[
A_\ell\subseteq B_{\ell,\delta}^{+},\qquad
\dim_{\mathrm M} A_2=\dim_{\mathrm M} A_3=\alpha>1-\varepsilon.
\]
The sets $A_\ell$ satisfy the multiplicative invariance hypotheses of \cite[Definition 1.5]{GMR} in bases $2^L$ and $3^L$.
\end{theorem}
\begin{proof}
Choose $m\ge1$ with $2^{-m}<\delta$. Let $L>m$ be a product of $r$ distinct odd primes, and let $Q_L$ be the square residues modulo $L$, including zero. Define residue classes
\[
F=\bigcup_{t=1}^{m}(Q_L-t),\qquad
J=\{0\le a<L:a\bmod L\notin F\}.
\]
By the Chinese remainder theorem,
\[
\frac{|Q_L|}{L}
=\prod_{p\mid L}\frac{p+1}{2p}\le(2/3)^r,
\qquad
\frac{|F|}{L}\le m(2/3)^r.
\]
Choose $r$ large enough that $|F|/L<\varepsilon$. Put $b_\ell=\ell^L$ and
\[
D_\ell=\left\{\sum_{a\in J}d_a\ell^a:
                     d_a\in\{0,\ldots,\ell-1\}\right\}.
\]
Every integer with base-$b_\ell$ digits in $D_\ell$ has zero base-$\ell$ digits at every nonnegative position $n^2-t$, $1\le t\le m$. If $n^2\ge m$, its residue modulo $\ell^{n^2}$ is therefore less than $\ell^{n^2-m}$; if $n^2<m$, that residue is zero. For each fixed integer the remaining coordinates tend to zero, so its one-prime norm is strictly less than $\ell^{-m}<\delta$.

The alphabets have cardinality $\ell^{|J|}$. Counting words of length $N$ gives
\[
|A_\ell\cap[0,b_\ell^N)|=|D_\ell|^N,
\qquad
\alpha=\frac{\log|D_\ell|}{\log b_\ell}=\frac{|J|}{L}>1-\varepsilon.
\]
Deleting either endpoint base-$b_\ell$ digit preserves the alphabet, which is precisely the invariance required in \cite{GMR}.
\end{proof}

\begin{corollary}[Consequences of transversality]\label{cor:gmr}
For the sets in Theorem~\ref{thm:stationary}, every $\zeta>0$ and all sufficiently large $X$ satisfy
\[
|A_2\cap A_3\cap[0,X)|\le X^{2\alpha-1+\zeta},
\qquad
\dim_{\mathrm M}(A_2+A_3)=1.
\]
\end{corollary}
\begin{proof}
The bases $2^L$ and $3^L$ are multiplicatively independent. Apply \cite[Theorems B and C]{GMR}; the mass and discrete Hausdorff dimensions of these restricted digit sets agree, and $\alpha>1/2$.
\end{proof}
These are upper bounds for intersections and lower bounds for sumsets. An intersection upper bound permits the intersection to consist only of zero; a sum of separate one-prime witnesses need not satisfy both constraints. Neither conclusion supplies an element of $B_\delta$ outside $p\Z$.

\begin{proposition}[Thickness obstruction]\label{prop:thin}
Let $b=\ell^L$, where $\ell\ge2$ and $L\ge1$ are integers. Let $D\subseteq\{0,\ldots,b-1\}$ be nonempty, and assume every finite $D$-word represents an integer in $B_{\ell,\delta}^{+}$, with the same square-exponent norm and $0<\delta<1/3$. Its associated real Cantor set
\[
K_D=\left\{\sum_{j\ge1}d_jb^{-j}:d_j\in D\right\}
\]
satisfies
\[
\tau(K_D)\le\frac{\delta}{1-2\delta}<1.
\]
Here $\tau$ is the thickness of \cite[Definition 8]{YY}, and a singleton has thickness zero.
\end{proposition}
\begin{proof}
For prefix lengths $N=Lt^2$ we have $b^N=\ell^{(Lt)^2}$. Apply the square-scale constraint to the integer represented by each prefix and take limits. This gives
\[
K_D\subseteq[0,\delta]\cup[1-\delta,1].
\]
If $K_D$ meets both bands, its central complementary gap has length at least $1-2\delta$. Every other gap lies inside a band and has length at most $\delta$, so the central gap is largest. The bridges on its two sides have lengths at most $\delta$. The definition of thickness gives the claimed bound.

If $K_D$ lies in one band, translate its minimum to zero. The translated alphabet is $D-\min D$, and the diameter $T$ of the translated Cantor set is at most $\delta$. If it is not a singleton, let $d\ge1$ be the smallest positive translated digit. The gap after the first cylinder has length $(d-T)/b$. Its left bridge has length at most $T/b$; no gap inside that cylinder is removed earlier, since its length is at most $T/b<(d-T)/b$. Hence
\[
\tau(K_D)\le\frac{T}{d-T}
\le\frac{\delta}{1-\delta}
\le\frac{\delta}{1-2\delta}.
\]
This includes every case.
\end{proof}
Thickness is invariant under affine changes of scale. Thus grouping digits into arbitrarily large power bases does not give high thickness for stationary inner models. In dimension one the hypothesis of \cite[Theorem 11]{YY} requires a sufficiently small sum of negative powers of the thicknesses, with constant $432^2$; thicknesses below one cannot satisfy it. This excludes that direct stationary-inner-model application, not every possible nonstationary or higher-dimensional adaptation.

A lower recurrence or intersection estimate with residue-class control remains missing. Theorem~\ref{thm:density} also shows why finitely many proper stationary inner models cannot cover the whole mixed ball: their mass dimensions are below one.

\section{Sparse series and single-scale saturation}\label{sec:sparse}
Throughout this section $S=\{2,3\}$.
Let $n_s$ be strictly increasing positive integers, and put
\[
q_s=2^{n_s^2},\qquad
x=\sum_{s\ge1}q_s^{-1},\qquad
Q_J=q_{J+2},\qquad
T_J=\sum_{s\ge J+2}q_s^{-1}.
\]
The ratio of successive denominators is at least $8$, so $Q_J^{-1}\le T_J\le2Q_J^{-1}$.

\begin{lemma}[Relative tail estimate]\label{lem:sparse-tail}
If $q_J\mid k$, then
\[
kx\equiv k/q_{J+1}+kT_J\pmod1,\qquad
\left|\nrm{kx}-\nrm{k/q_{J+1}}\right|\le2|k|/Q_J.
\]
For $1\le|k|\le K$,
\[
P(k)\ge2^{-2\sqrt{\log_2(2|k|)}-2},
\qquad
\frac{|k|T_J}{P(k)}
\le\frac{8K}{Q_J}2^{2\sqrt{\log_2(2K)}}.
\]
The last bound tends to zero uniformly for $K\le Q_J^\theta$ and any fixed $0<\theta<1$.
\end{lemma}
\begin{proof}
Divisibility clears every term through $J$, and the circle norm is $1$-Lipschitz. For the lower bound on $P(k)$, choose $n=\lceil\sqrt{\log_2(2|k|)}\rceil$. Then $2^{n^2}\ge2|k|$ and
\[
P(k)\ge |k|/2^{n^2}
\ge2^{-2\sqrt{\log_2(2|k|)}-2}.
\]
Combine this with $T_J\le2/Q_J$. When $K\le Q_J^\theta$, the resulting factor is a negative power of $Q_J$ times $2^{O(\sqrt{\log Q_J})}$, so tends to zero.
\end{proof}

\begin{proposition}[Conditional exclusion of sparse parameters]\label{prop:sparse-tests}
Suppose along infinitely many $J$ there are nonzero integers $k_J$ with $q_J\mid k_J$ and $P(k_J)\to0$. Either of the following additional conditions implies $x\notin H$.

(a) For fixed $0<\theta<1$ and $\kappa>0$,
\[
|k_J|\le Q_J^\theta,\qquad
\nrm{k_J/q_{J+1}}\ge\kappa P(k_J).
\]

(b) For fixed $\kappa>0$, the witnesses are positive and
\[
k_J\le Q_J/4,\qquad
\kappa P(k_J)\le\{k_J/q_{J+1}\}\le P(k_J).
\]
Here braces denote fractional part.
If $n_{s+1}/n_s\to\infty$, condition (a) permits the size bound
$|k_J|\le q_{J+1}^A$ for any fixed $A>0$ in place of $Q_J^\theta$.
\end{proposition}
\begin{proof}
In (a), Lemma~\ref{lem:sparse-tail} makes the tail error $o(P(k_J))$, contradicting the flatness criterion of Theorem~\ref{thm:flatness} if $x\in H$. Under the final growth assumption, $A n_{J+1}^2\le n_{J+2}^2/2$ eventually, so $\theta=1/2$ suffices.

In (b), the tested coordinate at $Q_J$ gives $k_J/Q_J\le P(k_J)$, hence $0\le k_JT_J\le2P(k_J)$. After the integer prefix has been removed, the real sum of the active fractional part and the tail lies between $\kappa P(k_J)$ and $3P(k_J)$. Eventually it is below $1/2$, so no cancellation or wraparound occurs. Its norm is at least $\kappa P(k_J)$, again contradicting flatness.
\end{proof}
The prefix condition cannot be replaced merely by $P(k_J)\to0$; a sufficient explicit replacement is $P(k_J)<2^{-(2n_J-1)}$, by Proposition~\ref{prop:divisibility-radii}. The larger size cutoff in (b), without its one-sided condition, controls the relative tail only by a constant and does not prevent cancellation. The existence of saturation witnesses satisfying (a) or (b) is not proved.

For $n_s=s!$, the first three terms have dyadic denominator $2^{36}$ and the next denominator is $2^{576}$. On $|k|\le10^{10}$, the character values differ from those of the dyadic truncation by at most $2\cdot10^{10}/2^{576}$. Such finite experiments cannot separate this irrational from that truncation.

\begin{remark}[Near-resonances]\label{rem:sparse-resonance}
For $x=\sum_{s\ge1}2^{-(s!)^2}$, let $K_j=2^{(j!)^2}$ and $\nu_j=-\sum_{s=1}^j2^{-(s!)^2}$. Then
\[
0<\nrm{x+\nu_j}\le2K_j^{-(j+1)^2}.
\]
The number is irrational, since its binary expansion is not eventually periodic. These permitted finite combinations of square-position frequencies show that a uniform lower bound of size $K^{-1+\eta}$ for every nontrivial tuple and every irrational parameter cannot hold. A Fourier approach must control such near-resonances rather than exclude them all.
\end{remark}

\section{Joint approximation and splitting of characters}\label{sec:joint}
In this section write $P_2=P_{\{2\}}$, $P_3=P_{\{3\}}$ and $P=\max(P_2,P_3)$, and let $H_2,H_3,H$ be their continuous parameter groups. The following arithmetic property is sufficient to determine $H$, but is not proved here.

\begin{definition}[Joint approximation]\label{def:JA}
For $t\in\Z$, let $\mathrm{JA}_t$ denote the assertion that, for every $\delta>0$, there is $k\in\Z$ with
\[
P_3(k)<\delta,\qquad P_2(k-t)<\delta.
\]
We write $\mathrm{JA}$ for $\mathrm{JA}_1$.
\end{definition}

\begin{proposition}[Equivalent formulations]\label{prop:JA-equivalent}
Put $U_\delta=\{u:P_2(u)<\delta\}$ and $V_\delta=\{v:P_3(v)<\delta\}$. The following are equivalent:
\begin{enumerate}
\item $\mathrm{JA}$ holds;
\item there are integers $k_j$ with $P_3(k_j)\to0$ and $P_2(k_j-1)\to0$;
\item $U_\delta-V_\delta=\Z$ for every $\delta>0$;
\item the diagonal integers are dense in $G_2\times G_3$, where $G_p$ is the metric completion of $(\Z,P_p)$.
\end{enumerate}
Under these conditions the completion of $(\Z,P)$ is $G_2\times G_3$.
\end{proposition}
\begin{proof}
The first two formulations are equivalent by taking radii tending to zero. Given the second, a target $t$ is realized by
$u=-t(k_j-1)$ and $v=-tk_j$; both norms tend to zero and $u-v=t$. Conversely, target $1$ in the third formulation gives $k=-v$ in the first.

Let $D$ be the closure of the diagonal integers in the product. The second formulation puts $(1,0)$ in $D$. Since $(1,1)\in D$, also $(0,1)\in D$, and therefore $\Z^2\subseteq D$. Integer pairs are dense in the product, proving the fourth formulation. The converse follows by approximating $(1,0)$. Finally, the diagonal embedding is isometric for the maximum product metric, so its closure is the mixed completion.
\end{proof}

\begin{theorem}[Conditional splitting]\label{thm:JA-split}
Assume $\mathrm{JA}$. If $x\in H$ satisfies $\nrm{mx}\le CP(m)$ for every $m\in\Z$, with $C\ge0$, then there is $y\in\T$ such that
\[
\nrm{my}\le CP_2(m),\qquad
\nrm{m(x-y)}\le CP_3(m)\qquad(m\in\Z).
\]
Consequently
\[
H=H_2+H_3=\Z[1/6]/\Z.
\]
In particular, $\mathrm{JA}$ excludes both irrational continuous parameters and rational denominators involving a prime outside $\{2,3\}$.
\end{theorem}
\begin{proof}
Choose $k_j$ as in Proposition~\ref{prop:JA-equivalent}. Compactness of $\T$ gives one subsequence along which $k_jx\to y$. Put $a_j=P_2(k_j-1)$ and $b_j=P_3(k_j)$. Subadditivity gives, for each integer $m$,
\[
\begin{aligned}
P(mk_j)&\le\max\{P_2(m)+|m|a_j,\ |m|b_j\},\\
P(m(1-k_j))&\le\max\{|m|a_j,\ P_3(m)+|m|b_j\}.
\end{aligned}
\]
Apply the Lipschitz bound to these two integers and pass to the limit along the chosen subsequence. Multiplication by each fixed $m$ and the norm on $\T$ are continuous, giving the stated inequalities. The subsequence does not depend on $m$. Thus $y\in H_2$ and $x-y\in H_3$. Theorem~\ref{thm:oneprime} gives the upper inclusion in the displayed equality; the reverse inclusion follows because the mixed topology makes both one-prime character groups continuous.
\end{proof}

\begin{lemma}[Approximation by multiples of a unit]\label{lem:unit-approx}
Let $p$ be a prime, $d\ge1$ with $\gcd(d,p)=1$, and $\eta>0$. There is an integer $a\ge1$ with $P_{\{p\}}(da-1)<\eta$.
\end{lemma}
\begin{proof}
For $d=1$ take $a=1$. Otherwise choose $N\ge1$ so $(d-1)/p^{2N+1}<\eta$, and choose $0\le w<d$ with $p^{N^2}w\equiv-1\pmod d$. Set $a=(p^{N^2}w+1)/d$. The square coordinates through $N$ vanish, while every later coordinate has norm at most its nonnegative real value. Therefore
\[
P_{\{p\}}(da-1)\le w/p^{2N+1}<\eta.
\]
\end{proof}

\begin{theorem}[Any nonzero target suffices]\label{thm:nonzero-target}
For every fixed $t\in\Z\setminus\{0\}$, $\mathrm{JA}_t$ is equivalent to $\mathrm{JA}$. It also suffices that, for every $\delta>0$, some target $t_\delta$ with $0<|t_\delta|\le T_0$ is attainable, for a fixed finite bound $T_0$.
\end{theorem}
\begin{proof}
Changing signs reduces to $t>0$. Write $t=ds$, where $d=2^{v_2(t)}$ and $s$ is odd. Fix $\delta>0$. By Lemma~\ref{lem:unit-approx} in base three, choose $a\ge1$ with $P_3(da-1)<\delta/(2s)$. Apply $\mathrm{JA}_t$ at precision $\delta/(2a)$ and put
\[
k'=ak+s(1-ad).
\]
Then $k'-s=a(k-t)$, so $P_2(k'-s)<\delta/2$, while
\[
P_3(k')\le aP_3(k)+sP_3(1-ad)<\delta.
\]
This proves $\mathrm{JA}_s$. Choose $b\ge1$ with $P_2(sb-1)<\delta/2$, and apply $\mathrm{JA}_s$ at precision $\delta/(2b)$ to obtain $h$. Then
\[
P_3(bh)<\delta/2,\qquad
P_2(bh-1)\le bP_2(h-s)+P_2(bs-1)<\delta.
\]
Thus $\mathrm{JA}$ holds. The converse follows by multiplying its witnesses by $t$ and using precision $\delta/|t|$. For the bounded family, pigeonhole gives a fixed nonzero target along radii tending to zero; monotonicity supplies its approximation at every radius.
\end{proof}

\begin{corollary}[Dichotomy of attainable shifts]\label{cor:shift-dichotomy}
The set $I=\bigcap_{\delta>0}(U_\delta-V_\delta)$ is either $\{0\}$ or $\Z$.
\end{corollary}
\begin{proof}
By symmetry, $t\in I$ means $\mathrm{JA}_t$. Zero always belongs to $I$, and Theorem~\ref{thm:nonzero-target} applies to any nonzero member.
\end{proof}
Theorem~\ref{thm:JA-split} also applies to a partition of a finite prime set into two nonempty parts, provided the analogous approximation between the two parts is assumed. Induction requires this property for the successive partitions; approximation only for individual pairs is not asserted to suffice. We do not prove that $\mathrm{JA}$ is necessary for the equality $H=\Z[1/6]/\Z$.

\section{Residue distribution and partial block estimates}\label{sec:block-residues}
The results here describe individual digit models and a restricted mixing range. They do not establish the joint approximation property.

\begin{proposition}[Residues in the stationary inner models]\label{prop:stationary-residues}
Let $A_\ell$ be a stationary inner model from Theorem~\ref{thm:stationary}, for $\ell=2$ or $3$. For every fixed integer $q\ge1$ coprime to $\ell$ and every residue $r$,
\[
\frac{|A_\ell\cap[0,X)\cap(r+q\Z)|}{|A_\ell\cap[0,X)|}\longrightarrow\frac1q.
\]
Each residue slice has the same mass dimension as $A_\ell$.
\end{proposition}
\begin{proof}
The alphabet in base $\ell^L$ contains $0$ and $\ell^a$ for every free position $a\in J$, and $J$ is nonempty. Hence its digit differences, together with $q$, have greatest common divisor one. Corollary A of \cite{SA} gives uniform distribution. A positive limiting proportion preserves the mass dimension.
\end{proof}

\begin{lemma}[Quantitative one-base residue estimate]\label{lem:free-residues}
Fix a prime $\ell$, integers $m\ge1$, $N\ge1$, and $q\ge2$ coprime to $\ell$. Let $E_{\ell,m}(N)$ consist of the integers below $\ell^N$ whose digits are zero at every position $n^2-t\in[0,N-1]$, for $n\ge1$ and $1\le t\le m$. Let $f_N$ be the number of remaining free positions. Then
\[
f_N\ge N-m\lfloor\sqrt{N+m}\rfloor,\quad
|E_{\ell,m}(N)|=\ell^{f_N},\quad
E_{\ell,m}(N)\subset B_{\ell,\ell^{-m}}^+.
\]
Put
\[
\rho_{\ell,q}=
\max_{1\le c<q}\left|\frac1\ell\sum_{d=0}^{\ell-1}
e^{2\pi i cd/q}\right|<1.
\]
For every residue $r$,
\[
\left|\frac{|\{k\in E_{\ell,m}(N):k\equiv r\pmod q\}|}{\ell^{f_N}}-\frac1q\right|
\le\frac{q-1}{q}\rho_{\ell,q}^{f_N}.
\]
\end{lemma}
\begin{proof}
There are at most $m\lfloor\sqrt{N+m}\rfloor$ forced positions. The last $m$ digits before each square cutoff are zero, or the entire residue is zero, proving the norm inclusion as in Theorem~\ref{thm:stationary}. For a uniformly chosen member, each nonzero Fourier coefficient modulo $q$ is the product over free positions $j$ of
$\ell^{-1}\sum_{d=0}^{\ell-1}e^{2\pi i hd\ell^j/q}$. Invertibility of $\ell$ modulo $q$ bounds every factor by $\rho_{\ell,q}$. The terms for digits zero and one differ, so $\rho_{\ell,q}<1$. Fourier inversion gives the estimate.
\end{proof}
The modulus is fixed here: as $q\to\infty$, $\rho_{\ell,q}\ge1-O_\ell(q^{-2})$. The estimate cannot simply be iterated with the other base's rapidly growing square denominators. Separate residue equidistribution does not determine intersections: the sets with residues $\{0,\ldots,p-1\}$ and $\{0,p+1,\ldots,2p-1\}$ modulo $2p$ are both uniformly distributed modulo $p$, but their intersection is $2p\Z$ among nonnegative integers.

\begin{proposition}[Obstruction for small initial alphabets]\label{prop:monomial}
For every prime $\ell$ and $L\ge1$, every stationary alphabet in base $b=\ell^L$ containing $0$ and $1$ contains infinitely many integers $k$ with $P_{\{\ell\}}(k)=1/\ell$.
\end{proposition}
\begin{proof}
For $t\ge0$ put $n=Lt+1$, $j=Lt^2+2t$, and
$k=b^j=\ell^{n^2-1}$. Its expansion has just one digit $1$. The square coordinates before $n$ vanish; at $n$ the norm is $1/\ell$, and the later values are smaller.
\end{proof}
Thus initial small-digit alphabets in arbitrarily large power bases are not inner models for radii below $1/\ell$. The simultaneous small-digit theorem and block method of \cite{BC} allow exceptional digits or blocks and do not directly impose every prescribed square window. This obstructs that substitution, not every adaptation of their method.

\begin{proposition}[Fourier contraction across free ternary blocks]\label{prop:block-fourier}
Let $F$ be a finite set of nonnegative digit positions and
$K=\sum_{j\in F}\varepsilon_j3^j$, with independent uniform digits $\varepsilon_j\in\{0,1,2\}$. Partition $F$ into maximal consecutive intervals $I$. For $1\le h<2^t$, write $h/2^t=a/2^s$ in lowest terms and set $L_s=\lceil s\log_3 2\rceil+1$. Then
\[
\left|\mathbb E e^{2\pi i hK/2^t}\right|
\le(2/3)^{\sum_I\lfloor |I|/L_s\rfloor}.
\]
If $f=|F|$ and there are $r$ intervals, the exponent is at least $f/L_s-r$.
\end{proposition}
\begin{proof}
Put $\theta_j=a3^j/2^s$ modulo one. Its norm is at least $2^{-s}$, and while below $1/6$ it triples at the next position. Consequently every $L_s$ consecutive positions include one with $\nrm{\theta_j}\ge1/6$. At that position the normalized digit factor satisfies
\[
\left|\frac{1+e^{2\pi i\theta_j}+e^{4\pi i\theta_j}}3\right|
=\frac{|1+2\cos(2\pi\theta_j)|}{3}\le2/3.
\]
All other factors have modulus at most one. Multiplication over disjoint subintervals of length $L_s$ proves the claim, and $\lfloor u\rfloor\ge u-1$ gives the last bound.
\end{proof}
For the ternary zero-window model through $M^2$ positions, $f=M^2-O_\delta(M)$ and $r=O_\delta(M)$. The proposition, together with Fourier inversion, gives useful uniform residue estimates for $t\le c_\delta M$ with a sufficiently small constant. It gives no such conclusion when $t$ is of order $M^2$. Uniformity at every highest modulus is stronger than required and can fail by cardinality: a distribution supported on at most $3^f$ residues has total variation distance at least $1-3^f/2^t$ from uniform measure if $3^f<2^t$. These models force zero windows, so constitute sufficient constructions rather than equivalent descriptions of the norm balls.

Short-orbit exponential-sum estimates such as \cite{Vandehey} provide relevant methods, but their applicability to the joint, conditioned square windows requires further estimates. The restricted bound above is not an impossibility theorem for Markov or block methods.

\section{Survival counts and the unresolved late scales}\label{sec:survival}
Here $S=\{2,3\}$, $0<\delta<1/2$, and $X\ge1$ is an integer. Define
\[
N_\delta(X)=|\{1\le k\le X:P_3(k)<\delta,\ P_2(k-1)<\delta\}|.
\]
Only denominators $p^{n^2}\le X/\delta$ require testing. At every later scale the relevant nonnegative numerator is at most $X$, so its real ratio is strictly below $\delta$. Order these active denominators as $q_1<\cdots<q_J$, remembering the prime of each constraint, and let $A_j(X)$ count candidates surviving the first $j$ tests; thus $A_0(X)=X$ and $A_J(X)=N_\delta(X)$. Put $Q_j=\operatorname{lcm}(q_1,\ldots,q_j)$, and let $R_j$ count allowed residues modulo $Q_j$, with $Q_0=R_0=1$.

\begin{proposition}[Complete periods and early survival]\label{prop:early-survival}
For every $j$,
\[
\lfloor X/Q_j\rfloor R_j\le A_j(X)\le\lceil X/Q_j\rceil R_j.
\]
If $\alpha_j$ is the ratio of the new denominator to the preceding denominator of the same prime (using $1$ before its first appearance), then
\[
2\delta-\frac2{\alpha_j}\le
\frac{R_j/Q_j}{R_{j-1}/Q_{j-1}}
\le2\delta+\frac2{\alpha_j}.
\]
After a fixed initial prefix depending on $\delta$, each step with $Q_j\le X$ satisfies
\[
A_j(X)\ge(\delta/4)A_{j-1}(X).
\]
At the last such step $j_e$ there is a constant $C_\delta$ with
\[
A_{j_e}(X)\ge X\exp(-C_\delta\sqrt{\log X})
\]
for all sufficiently large $X$.
\end{proposition}
\begin{proof}
Partition the interval into full periods and a final partial period. For each previously allowed residue, its $\alpha_j$ lifts have one-prime phases equally spaced by $1/\alpha_j$. An arc of length $2\delta$ admits between $2\delta\alpha_j-2$ and $2\delta\alpha_j+2$ lifts. CRT leaves the other prime component unrestricted, giving the density bounds. Every $R_j$ is positive: $k\equiv1$ in the binary component and $k\equiv0$ in the ternary component satisfy all finite tests.

Eventually $\alpha_j\ge2/\delta$, so the complete-period density ratio is at least $\delta$. If $Q_j\le X$, then also $Q_{j-1}\le X$, and the period bounds imply
\[
A_j(X)\ge\frac X2\frac{R_j}{Q_j},\qquad
A_{j-1}(X)\le2X\frac{R_{j-1}}{Q_{j-1}}.
\]
This proves the safe factor $\delta/4$. The finite initial prefix has positive density, and the number of steps is $O_\delta(\sqrt{\log X})$, proving the last assertion.
\end{proof}
An exact full-period ratio need not equal the truncated ratio. In particular $Q_j\le X$ alone does not remove the partial period.

\begin{proposition}[Location and number of late steps]\label{prop:late-scales}
Put $a=(\log2)^{-1/2}+(\log3)^{-1/2}$, with natural logarithms. Then
\[
J=a\sqrt{\log X}+O_\delta(1),\qquad
J-j_e=(1-1/\sqrt2)a\sqrt{\log X}+O_\delta(1).
\]
The transition $Q_j\approx X$ occurs at an individual denominator
$q_j=X^{1/2}\exp(O(\sqrt{\log X}))$. Thus the late steps form asymptotically $1-1/\sqrt2$, about $29.3\%$, of all steps.
\end{proposition}
\begin{proof}
The number of square denominators at most $T$ is
$a\sqrt{\log T}+O(1)$. At a current denominator $q$, the largest denominator in each prime chain has logarithm $\log q+O(\sqrt{\log q})$, by taking the integer part of $\sqrt{\log_p q}$. Hence $\log Q_j=2\log q+O(\sqrt{\log q})$. Consecutive changes in $\log Q_j$ are also $O(\sqrt{\log q})$. At its crossing of $\log X$, therefore, $\log q=\tfrac12\log X+O(\sqrt{\log X})$. Substitution into the counting formula gives $j_e=a\sqrt{\log X}/\sqrt2+O_\delta(1)$ and the result. The active upper cutoff is $X/\delta$, which changes the count by $O_\delta(1)$.
\end{proof}
The late regime is not confined to denominators between $X$ and $X/\delta$; that interval contains only $O_\delta(1)$ scales.

\begin{theorem}[A sufficient survival estimate]\label{thm:survival}
Fix $\delta\in(0,1/2)$. Suppose constants $c_\delta\in(0,1)$ and $B_\delta\ge0$, independent of $X$ and $j$, satisfy
\[
A_j(X)\ge c_\delta A_{j-1}(X)-B_\delta
\]
for every active step and all sufficiently large integer $X$. Then
\[
N_\delta(X)\ge c_\delta^J X-\frac{B_\delta}{1-c_\delta}
\ge X\exp(-C_\delta\sqrt{\log X})-O_\delta(1).
\]
It suffices to assume the recurrence only for late steps, using Proposition~\ref{prop:early-survival} as the starting bound, or for bounded groups of consecutive constraints. If such bounds hold for every $\delta\in(0,1/2)$, then $\mathrm{JA}$ holds and $H=\Z[1/6]/\Z$.
\end{theorem}
\begin{proof}
Iteration gives
$A_J\ge c_\delta^J X-B_\delta\sum_{i=0}^{J-1}c_\delta^i$.
Use the geometric sum and $J=O_\delta(\sqrt{\log X})$. The same calculation starts at $j_e$ with its established lower bound, or applies to groups of bounded length. The resulting bound tends to infinity. Thus a witness exists at every radius, and Theorem~\ref{thm:JA-split} applies.
\end{proof}
The late-step hypothesis is unproved. The theorem requires only controlled lower losses, not independence or a ratio near $2\delta$. It is a sufficient quantitative condition, not a condition shown equivalent to $\mathrm{JA}$.

\begin{remark}[Exact finite survival counts]\label{rem:survival-counts}
For $X=2\,000\,000\,000$, the active denominators at both $\delta=\frac{1}{10}$ and $\delta=\frac{1}{20}$ are those in the following table. Exact integer comparisons give the displayed counts.
\begin{center}
\begin{tabular}{rrr}
$q_j$ & $A_j(X)$, $\delta=\frac{1}{10}$ & $A_j(X)$, $\delta=\frac{1}{20}$\\
\hline
$2$ & $1\,000\,000\,000$ & $1\,000\,000\,000$\\
$3$ & $333\,333\,333$ & $333\,333\,333$\\
$16$ & $41\,666\,666$ & $41\,666\,666$\\
$81$ & $7\,716\,049$ & $4\,629\,629$\\
$512$ & $1\,687\,880$ & $434\,025$\\
$19\,683$ & $340\,359$ & $44\,653$\\
$65\,536$ & $66\,472$ & $4\,538$\\
$33\,554\,432$ & $13\,351$ & $452$\\
$43\,046\,721$ & $2\,733$ & $49$\\
\end{tabular}
\end{center}
The first three constraints reduce the candidates exactly to $k\equiv33\pmod{48}$. Enumerating that progression and testing $d\min(r,q-r)<q$, with $d=10$ or $20$, verifies the remaining rows without floating-point decisions. The next square denominator in each chain exceeds $X/\delta$, certifying the tails. Here $Q_7=1\,289\,945\,088<X$, but $Q_8=660\,451\,885\,056>X$.
These finite counts support the proposed lower bound, but do not establish uniform survival or asymptotic independence.
\end{remark}

\section{Other prime pairs and finite prime sets}\label{sec:general-primes}
The conditional results of Section~\ref{sec:joint} do not depend on the special choice of the primes $2$ and $3$. For distinct primes $p,q$, put
\[
F_{p,q}(k)=\max\{P_{\{q\}}(k),P_{\{p\}}(k-1)\},
\qquad J_{p,q}=\inf_{k\in\Z}F_{p,q}(k).
\]
Write $\mathrm{JA}^{p,q}_t$ if, for every $\delta>0$, there is $k\in\Z$ with $P_{\{q\}}(k)<\delta$ and $P_{\{p\}}(k-t)<\delta$.

\begin{theorem}[General prime-pair reduction]\label{thm:general-pair}
For every pair of distinct primes $p,q$:
\begin{enumerate}
\item $\mathrm{JA}^{p,q}_t$ for any fixed nonzero integer $t$ is equivalent to $\mathrm{JA}^{p,q}_1$;
\item $\mathrm{JA}^{p,q}_1$ implies $H_{\{p,q\}}=\Z[1/(pq)]/\Z$;
\item $\mathrm{JA}^{p,q}_1$ is equivalent to $J_{p,q}=0$, and $J_{p,q}=J_{q,p}\le1/\max(p,q)$.
\end{enumerate}
The bounded-nonzero-target criterion and attainable-shift dichotomy also hold.
\end{theorem}
\begin{proof}
For $t>0$ write $t=ds$, with $d=p^{v_p(t)}$ and $\gcd(s,p)=1$. The integer $d$ is coprime to $q$. Given $\delta>0$, Lemma~\ref{lem:unit-approx} in base $q$ supplies $a\ge1$ with $P_{\{q\}}(da-1)<\delta/(2s)$. Use $\mathrm{JA}^{p,q}_t$ at precision $\delta/(2a)$ and set $k'=ak+s(1-ad)$. Then
\[
P_{\{p\}}(k'-s)<\delta/2,\qquad
P_{\{q\}}(k')\le aP_{\{q\}}(k)+sP_{\{q\}}(1-ad)<\delta.
\]
This gives $\mathrm{JA}^{p,q}_s$. Lemma~\ref{lem:unit-approx} in base $p$, applied to $s$, completes the reduction to target one exactly as in Theorem~\ref{thm:nonzero-target}. Negative targets are treated by changing signs; multiplication by $t$ gives the reverse implication. Pigeonhole and monotonicity prove the bounded-target assertion.

The compactness argument of Theorem~\ref{thm:JA-split} splits every mixed character into one continuous parameter from each prime. Apply Theorem~\ref{thm:oneprime}. The equivalence with the infimum is immediate from its definition. Finally,
$F_{p,q}(k)=F_{q,p}(1-k)$ proves symmetry, while $k=0,1$ give the upper bound.
\end{proof}
Adding primes gives the inclusion of parameter groups in Corollary~\ref{cor:twoprimes}. Replacing a prime by a larger one does not give a proved comparison of the topologies. Fewer denominators below a height bound is a counting observation, not a topology inclusion. Also, positive finite-range searches are not symmetric under exchanging the two primes: $1\le k\le X$ becomes $1-X\le1-k\le0$.

\begin{theorem}[Early survival and late scales for finitely many primes]\label{thm:general-counting}
Let $S$ have $r\ge1$ distinct primes. At each prime choose a fixed target $a_\ell\in\{0,1\}$, and test
$\nrm{(k-a_\ell)/\ell^{n^2}}<\delta$ on $1\le k\le X$, where $0<\delta<1/2$. Order the active denominators $\ell^{n^2}\le X/\delta$, and define $A_j,Q_j,R_j$ as in Section~\ref{sec:survival}. Put
\[
a_S=\sum_{\ell\in S}(\log\ell)^{-1/2}.
\]
After a fixed initial prefix, every step with $Q_j\le X$ satisfies $A_j(X)\ge(\delta/4)A_{j-1}(X)$. If $j_e$ is the last such step and $J$ is the total number of active tests, then
\[
\begin{aligned}
J&=a_S\sqrt{\log X}+O_{S,\delta}(1),\\
j_e&=\frac{a_S}{\sqrt r}\sqrt{\log X}+O_{S,\delta}(1),\\
J-j_e&=(1-1/\sqrt r)a_S\sqrt{\log X}+O_{S,\delta}(1).
\end{aligned}
\]
The crossing $Q_j\approx X$ occurs at an individual denominator
$X^{1/r}\exp(O_S(\sqrt{\log X}))$.
\end{theorem}
\begin{proof}
Each new scale changes just one prime component, with same-prime ratio $\ell^{2n-1}$. CRT leaves the remaining components independent in a complete period. The density ratio is therefore between $2\delta-2/\ell^{2n-1}$ and $2\delta+2/\ell^{2n-1}$. After a finite prefix it is at least $\delta$. Every finite list is consistent, by choosing $k\equiv a_\ell$ in each prime-power component. The floor and ceiling estimates of Proposition~\ref{prop:early-survival} now give the same factor $\delta/4$.

At an individual denominator $D$, the number of scales so far is $a_S\sqrt{\log D}+O_S(1)$. The largest scale of each prime has logarithm $\log D+O_S(\sqrt{\log D})$, hence
\[
\log Q_j=r\log D+O_S(\sqrt{\log D}).
\]
The change at each step is also $O_S(\sqrt{\log D})$. Crossing $\log X$ therefore gives
$\log D=(\log X)/r+O_S(\sqrt{\log X})$. Substitute in the scale count. The upper cutoff $X/\delta$ changes the count only by $O_{S,\delta}(1)$.
\end{proof}
Thus the asymptotic fraction of late tests is $1-1/\sqrt r$, approximately $29.3\%$ for two primes and $42.3\%$ for three. This says nothing by itself about survival through those tests. Theorem~\ref{thm:general-pair} does not supply joint approximation for a partition into mixed blocks of primes; the single-prime unit lemma cannot simply be assumed for such a block.

\section{Certified interval construction and further witnesses}\label{sec:interval-search}
Return to $S=\{2,3\}$ and write $F(k)=F_{2,3}(k)$. For an integer $D\ge2$, let $K_D$ be the least positive integer with $F(k)<1/D$, if one exists, and put $K_D=+\infty$ otherwise. The following exact arithmetic rule replaces long scans by certified exclusions of intervals.

\begin{lemma}[Forbidden-interval jump]\label{lem:gap-jump}
For a coordinate with denominator $q$ and shift $a\in\{0,1\}$ put
$T=\lfloor(q-1)/D\rfloor$. Its strict radius-$1/D$ condition is equivalent to
\[
s=(k-a)\bmod q\in[0,T]\cup[q-T,q-1].
\]
If $T<s<q-T$, put $k'=k+q-T-s$. Every integer $u$ with $k\le u<k'$ fails this coordinate, and $k'$ passes it. Among several failed coordinates one may take the largest proposed jump without skipping a possible full witness.
\end{lemma}
\begin{proof}
The condition is $D\min(s,q-s)<q$, equivalently $\min(s,q-s)\le T$. In the forbidden interval the residues from $k$ through $k'-1$ increase from $s$ to $q-T-1$ without wrapping; the residue at $k'$ is $q-T$ modulo $q$ (zero when $T=0$). For the largest proposed jump, its own failed coordinate excludes the whole skipped interval.
\end{proof}

A finite prefix of coordinates has a complete finite set of admissible residues modulo their least common multiple. CRT constructs this list; jumping to its next member excludes only integers failing a prefix coordinate.

\begin{proposition}[Correctness of interval propagation]\label{prop:gap-search}
Start at $k=1$. Repeatedly move to the next admissible prefix residue, test the remaining square denominators at most $Dk$, and, if any fail, apply the largest jump from Lemma~\ref{lem:gap-jump}. If none fail, return $k$. Every update excludes all intervening positive integers. A returned value equals $K_D$, and if $K_D<+\infty$, the procedure terminates at it.
\end{proposition}
\begin{proof}
Every prefix jump and every gap jump preserves the assertion that all smaller positive integers have been excluded. A failed-coordinate jump is strictly increasing. Denominators greater than $Dk$ are automatic, since both relevant numerators are nonnegative and at most $k$. Thus a returned candidate is a full witness and is least. If a least witness exists, no update can skip it, and a strictly increasing integer sequence below it has only finitely many steps.
\end{proof}
This is a certified construction procedure, not a proof of termination for every $D$. Joint approximation may be tested using positive witnesses without strengthening the conjecture: each mixed ball contains arbitrarily large positive integers, since it contains nonzero integers at arbitrarily small radii and every finite nonzero set has a positive minimum norm. Given any integer witness at half the desired radius, add a sufficiently large positive mixed-small integer at the other half-radius.

For two selected top scales $A=2^{N^2}$ and $B=3^{M^2}$, put $T_A=\lfloor(A-1)/D\rfloor$ and $T_B=\lfloor(B-1)/D\rfloor$. The integer bands centered at $1+uA$ and $vB$ intersect exactly when
\[
|1+uA-vB|\le T_A+T_B.
\]
Their endpoints are integers, so ordinary interval overlap is enough. This exposes an inhomogeneous Diophantine approximation problem for possible band centers. It does not impose the lower square-position constraints, which still have to be checked.

\begin{proposition}[Computational certificates]\label{prop:new-witnesses}
Exact interval propagation and independent certificate checks give
\[
K_{81}=K_{100}=37074752603481701889.
\]
For this integer,
\[
\begin{aligned}
P_2(k-1)&=\frac{354032140747263}{36028797018963968},\\
P_3(k)&=\frac{17008458307642}{1853020188851841}.
\end{aligned}
\]
Both values are below $\frac{1}{100}$. Further,
\[
K_{200}=24168444555598036601864193,
\]
with
\[
P_2(k-1)=\frac{1305689979903}{281474976710656},
\qquad P_3(k)=\frac{8}{2187}.
\]
Finally,
\[
K_{500}=6207496413615241618428036085462998253569,
\]
with
\[
P_2(k-1)=\frac{1}{512},\qquad
P_3(k)=\frac{15152367064}{7625597484987}.
\]
In particular,
\[
J_{2,3}\le\frac{15152367064}{7625597484987}<\frac{1}{500}.
\]
\end{proposition}
\begin{proof}[Verification]
The first two certificates use all good binary residues modulo $2^{25}$ for $k-1$ and all good ternary residues modulo $3^{16}$ for $k$. Their CRT modulus is $1\,444\,408\,272\,617\,472$. At $D=81,100,200$ the seed lists have respectively $10\,491$, $7\,095$ and $315$ classes. Replaying the prefix and forbidden-interval updates certifies all smaller positive integers as excluded. The first two returned integers were also reproduced by an independent maximal-jump implementation with no CRT prefix. At $D=500$, that implementation returns the stated integer after $10\,604\,623$ iterations; the gap rule certifies minimality.

All displayed norms are independently computed with exact integer and rational arithmetic. For each numerator, the calculation includes the first square denominator exceeding twice its absolute value; beyond it the ratios are positive, decreasing and below $\frac12$, proving the tail bound. The largest coordinates of the first integer occur at binary index $8$ and ternary index $6$. For the second they occur at binary index $8$ and ternary index $4$, and for the third at binary index $5$ and ternary index $6$.
\end{proof}

\begin{proposition}[A symmetric finite minimum]\label{prop:symmetric-minimum}
One has the exact finite-range identity
\[
\min\{F(k):k\in\Z,\ -10^{18}<k\le10^{18}\}=\frac{1}{81}.
\]
\end{proposition}
\begin{proof}[Verification]
At strict radius $\frac{1}{81}$ the prefix sets have $39$ binary residues modulo $2^{25}$ and $269$ ternary residues modulo $3^{16}$, giving $10\,491$ CRT classes. They yield $7\,263\,183$ candidates in $1\le k\le10^{18}$. None survives the remaining active denominators
\[
2^{36},\quad3^{25},\quad2^{49},\quad3^{36},\quad2^{64}.
\]
All larger square denominators exceed $81\cdot10^{18}$ and are automatic. Reversing the shifts, so the binary test is on $k$ and the ternary test on $k-1$, gives $7\,263\,180$ positive candidates and again no survivor. Since
$F_{2,3}(k)=F_{3,2}(1-k)$, this covers $1-10^{18}\le k\le0$ as well. Exact integer enumeration was independently repeated in both orientations.

The integer $k=149176171345412097$ attains the bound:
\[
P_2(k-1)=\frac{84935169}{8589934592},\qquad
P_3(k)=\frac{1}{81}.
\]
\end{proof}
There is no conflict with Proposition~\ref{prop:new-witnesses}: its better witnesses lie beyond this finite range. A finite minimum is not a lower bound for the global infimum $J_{2,3}$.

Bounded interval runs also certify $K_{1000}>2\cdot10^{40}$; no value of $K_{1000}$ is asserted here. The absence of a witness below a finite cutoff cannot prove nonexistence. In particular, extrapolating a growth law from the listed values of $K_D$ would be a heuristic, not a proof of joint approximation.

\section{The quotient and compactness criteria}\label{sec:quotient}
Throughout this section $S=\{2,3\}$. Let $G_p$ be the metric completion of $(\Z,P_p)$ for $p=2,3$, and let $N_p$ be the continuous extension of $P_p$ to $G_p$. Put
\[
X=G_2\times G_3,\qquad N(x,y)=\max\{N_2(x),N_3(y)\}.
\]
Let $D$ be the closure in $X$ of the diagonal integers, and define
\[
Q=X/D,\qquad e=(1,0)+D,\qquad
\rho(z+D)=\inf_{d\in D}N(z-d).
\]
The metric of $Q$ is $d_Q(q,q')=\rho(q-q')$. For integer targets put
\[
J_m=\inf_{k\in\Z}\max\{P_2(k-m),P_3(k)\}.
\]
Thus $J_1=J_{2,3}$ in Section~\ref{sec:general-primes}. This quotient $Q$ is distinct from the circle parameter group $H$.

\begin{theorem}[The quotient obstruction]\label{thm:quotient-obstruction}
The group $Q$ is a Polish abelian group with a complete invariant metric, the cyclic group generated by $e$ is dense, and
\[
\rho(me)=J_m\qquad(m\in\Z),\qquad \widehat{Q}=\{0\}.
\]
The following statements are equivalent:
\begin{enumerate}
\item $\mathrm{JA}$ holds;
\item $Q=\{0\}$;
\item $Q$ is compact;
\item $Q$ is locally compact;
\item $Q$ is sigma-compact.
\end{enumerate}
\end{theorem}
\begin{proof}
The functions $P_p$ are Lipschitz for their own metrics and therefore extend to the norms $N_p$. Closedness of $D$ makes the quotient metric Hausdorff and gives the quotient topology: the inverse image of a radius-$r$ ball is $D+\{z:N(z)<r\}$. The quotient is separable. It is complete as well. Given a Cauchy sequence of cosets, pass to a subsequence whose successive distances are less than $2^{-j}$. Choose representatives inductively, adjusting by elements of $D$, so that successive representative distances are less than $2^{1-j}$. The representatives are Cauchy in the complete group $X$ and converge. The original Cauchy sequence of cosets then converges to the same coset limit. Hence $Q$ is Polish.

Every integer pair $(a,b)$ has the same coset as $(a-b,0)$, because $(b,b)\in D$. Integer pairs are dense in $X$, proving the density of the cyclic group generated by $e$. Since the diagonal integers are dense in $D$ and $N$ is continuous,
\[
\rho(me)=\inf_{k\in\Z}N((m,0)-(k,k))=J_m.
\]

A continuous character of $Q$ pulls back to a continuous character of $X$ vanishing on $D$. Every product character is the sum of its restrictions to the two factors. By Theorem~\ref{thm:oneprime}, the restrictions to the dense integer subgroups have parameters
\[
\alpha_2\in\Z[1/2]/\Z,\qquad
\alpha_3\in\Z[1/3]/\Z.
\]
Vanishing on every diagonal integer gives $\alpha_2+\alpha_3=0$ in $\T$. The two displayed subgroups intersect only at zero, so both restrictions vanish on the integers, hence on their completions. The pulled-back character and the original character are zero.

Proposition~\ref{prop:JA-equivalent} gives (1) if and only if (2), and (2) implies (3), (4), and (5). A locally compact Hausdorff abelian group has continuous characters separating its points, by Pontryagin duality \cite[Chapter 6]{Morris}. Since $\widehat{Q}=\{0\}$, (4) implies (2). Compactness implies local compactness. Finally, if $Q$ is a countable union of compact sets, Baire's theorem gives a compact set with nonempty interior. Translating that set gives a compact neighbourhood of zero, proving (5) implies (4).
\end{proof}

For $r>0$ define the attainable target differences
\[
U_{p,r}=\{k\in\Z:P_p(k)<r\},\qquad
W_r=U_{3,r}-U_{2,r}.
\]
Directly from the definition of the infimum,
\[
W_r=\{m\in\Z:J_m<r\}.
\]
A subset $A$ of $\Z$ is called \emph{syndetic} if $\Z=F+A$ for some finite $F\subset\Z$; equivalently, its gaps are bounded.

\begin{corollary}[A bounded-gap criterion]\label{cor:syndetic-JA}
Joint approximation holds if and only if $W_\varepsilon$ is syndetic for every $\varepsilon>0$. Explicitly,
\[
\mathrm{JA}\quad\Longleftrightarrow\quad
(\forall\varepsilon>0)(\exists F\subset\Z\ \text{finite})\
\Z=F+W_\varepsilon.
\]
\end{corollary}
\begin{proof}
If $\mathrm{JA}$ holds, Theorem~\ref{thm:nonzero-target} gives $J_m=0$ for every integer $m$, so $W_\varepsilon=\Z$ and $F=\{0\}$ suffices. Conversely, the displayed condition gives finite $\varepsilon$-nets for the dense cyclic subgroup of $Q$: if $m=f+w$ with $w\in W_\varepsilon$, then
\[
d_Q(me,fe)=J_{m-f}<\varepsilon.
\]
Given any desired radius $\eta>0$, apply this with $\varepsilon=\eta/2$ and approximate arbitrary points of $Q$ within $\eta/2$ by the dense cyclic subgroup. Thus $Q$ is totally bounded. Completeness makes $Q$ compact, and Theorem~\ref{thm:quotient-obstruction} gives $\mathrm{JA}$.
\end{proof}
The bound on the gaps may depend arbitrarily on $\varepsilon$. Syndeticity at a single fixed precision is not the assertion of Corollary~\ref{cor:syndetic-JA}.

\begin{proposition}[A local finite-cover criterion]\label{prop:local-cover}
Joint approximation is equivalent to the following property:
\[
(\exists r>0)(\forall\varepsilon>0)
(\exists F\subset\Z\ \text{finite})\
W_r\subseteq F+W_\varepsilon.
\]
\end{proposition}
\begin{proof}
Under $\mathrm{JA}$ every $W_r$ is $\Z$, so the property holds. Conversely, fix such an $r$. The orbit points $me$ with $m\in W_r$ are dense in the open ball $\{q:\rho(q)<r\}$: approximate a point strictly inside the ball by the dense cyclic subgroup. The finite-cover property gives finite nets for those orbit points at every precision. Using half the desired radius, as in Corollary~\ref{cor:syndetic-JA}, shows that the open ball is totally bounded. Its closure is totally bounded and complete, hence compact. It contains a neighbourhood of zero, so $Q$ is locally compact. Theorem~\ref{thm:quotient-obstruction} gives $\mathrm{JA}$.
\end{proof}
Here only the outer radius $r$ is fixed; the finite covers are still required at every finer precision.

\begin{corollary}[Separated targets if approximation fails]\label{cor:separated-targets}
If $\mathrm{JA}$ fails, then for every $r>0$ there exist $\varepsilon>0$ and distinct integers $m_1,m_2,\ldots$ such that
\[
J_{m_i}<r,\qquad
J_{m_i-m_j}\ge\varepsilon\quad(i\ne j).
\]
Moreover, every compact subset of $Q$ is nowhere dense, and $Q$ is not sigma-compact.
\end{corollary}
\begin{proof}
Negate the finite-cover condition in Proposition~\ref{prop:local-cover}. For each fixed $r$ there is $\varepsilon>0$ such that no finite collection of translates of $W_\varepsilon$ covers $W_r$. Choose $m_i$ successively outside the translates centred at the previously chosen integers. Symmetry of $J_m$ gives the stated separation. If a compact subset of $Q$ had nonempty interior, translating it would give a compact neighbourhood and hence local compactness, contradicting Theorem~\ref{thm:quotient-obstruction}. Compact subsets are closed, so they are nowhere dense. The same theorem excludes sigma-compactness.
\end{proof}

\begin{proposition}[Every finite-coordinate obstruction vanishes]\label{prop:finite-quotient}
For positive integers $N,M$ and every $m\in\Z$, put
\[
J_m^{N,M}=
\inf_{k\in\Z}
\max\left\{
\max_{1\le n\le N}\nrm{(k-m)/2^{n^2}},
\max_{1\le n\le M}\nrm{k/3^{n^2}}
\right\}.
\]
Then $J_m^{N,M}=0$.
\end{proposition}
\begin{proof}
The Chinese remainder theorem supplies an integer $k$ with
\[
k\equiv m\pmod{2^{N^2}},\qquad
k\equiv0\pmod{3^{M^2}}.
\]
Every displayed coordinate is then zero.
\end{proof}
The integers supplied by this proposition depend on the cutoffs and may escape to infinity. It does not justify interchanging the infimum over $k$ with the infinite coordinate supremum. In particular, triviality of every finite-coordinate quotient does not establish triviality or compactness of $Q$. A compactness proof by this route still needs uniform control of the omitted coordinates, or a genuinely compact family of suitable representatives. The known sigma-compactness of $H$ also gives no sigma-compactness assertion for the different group $Q$.

\begin{remark}[Finite shifted-target certificates]\label{rem:shift-certificates}
Exact computations give, for each of the $24$ nonzero shifts $m\in\{-12,\ldots,12\}$ and each $D\in\{20,50,100\}$, a positive integer $k$ with
\[
P_2(k-m)<\frac{1}{D},\qquad P_3(k)<\frac{1}{D}.
\]
All $72$ returned witnesses were independently checked with integer arithmetic at every coordinate through $n=40$. For each numerator $x=k-m$ or $x=k$, the first omitted denominator $p^{41^2}$ exceeds $D|x|$; every subsequent coordinate is therefore strictly below $1/D$. The decimal digit ranges of the returned witnesses are $1$--$9$, $6$--$18$, and $16$--$21$ for $D=20,50,100$, respectively. In particular,
\[
\{-12,\ldots,12\}\subseteq W_{1/100},
\]
with zero included trivially. These are finite full-norm certificates, not a proof of syndeticity: arbitrarily large gaps at the same precision could occur outside the observed window. No minimality claim for all $72$ integers is needed here.
\end{remark}

\section{Density and common Lipschitz criteria}\label{sec:density-JA}
In this and the next two sections the primes are $2$ and $3$, and $P_p(k)=\sup_{n\ge1}\nrm{k/p^{n^2}}$. We retain $U_{p,r}$, $W_r$, and $J_m$ from Section~\ref{sec:quotient}. For $A\subseteq\Z$ define its upper Banach density by
\[
d^*(A)=\limsup_{L\to\infty}\ \sup_{a\in\Z}
\frac{|A\cap\{a,\ldots,a+L-1\}|}{L}.
\]

\begin{lemma}[Positive density and difference sets]\label{lem:BD-difference}
If $d^*(A)=\alpha>0$, then $A-A$ is syndetic. More precisely, there is a finite $F$ with $\Z=F+(A-A)$ and $|F|\le1/\alpha$.
\end{lemma}
\begin{proof}
This is the standard difference-set lemma; see \cite{Beiglbock}. For completeness, take pairwise disjoint translates $f_1+A,\ldots,f_h+A$. On intervals of length tending to infinity where the density of $A$ tends to $\alpha$, each fixed translate has the same density up to a boundary error tending to zero. Disjointness gives $h\alpha\le1$. Greedily adjoining disjoint translates must therefore stop after at most $1/\alpha$ steps. For a maximal family, every $z+A$ meets one of its translates, so $z\in f_i+(A-A)$ for some $i$.
\end{proof}

\begin{theorem}[A density criterion for joint approximation]\label{thm:BD-JA}
Joint approximation holds if and only if $d^*(W_r)>0$ for every $r>0$. It suffices to verify this for a sequence of positive radii tending to zero. Equivalently, failure of joint approximation is equivalent to $d^*(W_r)=0$ for at least one positive radius.
\end{theorem}
\begin{proof}
Subadditivity gives $W_r-W_r\subseteq W_{2r}$. If $d^*(W_r)>0$, Lemma~\ref{lem:BD-difference} makes $W_{2r}$ syndetic. Applying this at arbitrarily small radii and then Corollary~\ref{cor:syndetic-JA} proves joint approximation. Conversely, joint approximation gives $W_r=\Z$ for all $r>0$. The last assertion is the logical negation of the equivalence.
\end{proof}
No uniform positive lower bound on the densities is required. The individual summands, however, have zero upper Banach density at small radii.

\begin{proposition}\label{prop:oneprime-BD}
For $p\in\{2,3\}$ and $0<r<\frac{1}{2}$, one has $d^*(U_{p,r})=0$.
\end{proposition}
\begin{proof}
Put $q_n=p^{n^2}$ and let $S_n$ be the good residues modulo $q_n$ for the first $n$ tests. Write $\beta_n=|S_n|/q_n$ and $\alpha_n=q_n/q_{n-1}$ for $n\ge2$. A fixed residue modulo $q_{n-1}$ has at most $2r\alpha_n+2$ lifts whose distance to zero modulo $q_n$ is less than $rq_n$: the admissible set is the union of two end intervals of total length $2rq_n$. Consequently
\[
\beta_n\le\beta_{n-1}\left(2r+\frac{2}{\alpha_n}\right).
\]
Since $\alpha_n\to\infty$, these factors are eventually bounded above by a constant less than one, and $\beta_n\to0$. The full ball is contained in each periodic prefix set, whose upper Banach density is $\beta_n$. This proves the claim.
\end{proof}
Thus a sumset theorem requiring positive density of the two individual balls does not directly give Theorem~\ref{thm:BD-JA}.

\begin{theorem}[A common Lipschitz certificate]\label{thm:common-Lipschitz}
Joint approximation fails if and only if there is a function $b:\Z\to\R$ with $b(1)\ne b(0)$ and finite constants $C_2,C_3\ge0$ such that, for all $m,n\in\Z$,
\[
|b(m)-b(n)|\le C_2P_2(m-n),\qquad
|b(m)-b(n)|\le C_3P_3(m-n).
\]
Every such function is bounded, and
\[
J_1\ge\frac{|b(1)-b(0)|}{C_2+C_3}>0.
\]
If joint approximation fails, the choice $b(m)=J_m$ works with $C_2=C_3=1$.
\end{theorem}
\begin{proof}
Apply the binary inequality between $1$ and $k$ and the ternary inequality between $k$ and $0$. Their sum bounds $|b(1)-b(0)|$ by $(C_2+C_3)\max\{P_2(k-1),P_3(k)\}$. Taking the infimum proves the lower bound. Boundedness follows from $P_p\le\frac{1}{2}$. Conversely the quotient norm gives
\[
|J_m-J_n|\le J_{m-n}\le\min\{P_2(m-n),P_3(m-n)\}.
\]
Here the last inequalities also follow by using $k=0$ and $k=m-n$ in the defining infimum. Failure of joint approximation means $J_1>0$, whereas $J_0=0$.
\end{proof}
This is an exact obstruction criterion, not the construction of an obstruction. The function need not be additive; the trivial continuous dual of the quotient does not force it to be constant.

\begin{proposition}[Two exclusions for candidate certificates]\label{prop:certificate-exclusions}
A function uniformly approximable both by functions of period a power of $2$ and by functions of period a power of $3$ is constant. Also, a function Lipschitz for one $P_p$ is constant if $b(n+1)-b(n)\to0$ as $n\to+\infty$.
\end{proposition}
\begin{proof}
For the first assertion choose periodic approximations $u,v$ with uniform errors at most $\eta$. Given $m,n$, the Chinese remainder theorem supplies $k$ congruent to $m$ modulo the period of $u$ and to $n$ modulo the period of $v$. Thus $|b(m)-b(n)|\le4\eta$. Let $\eta\downarrow0$.

For the second assertion, nondiscreteness of the one-prime metric gives positive integers $k_j$ with $P_p(k_j)\to0$ and $k_j\to\infty$ in the usual order. Indeed, choose arbitrarily small nonzero points, change their signs if needed, and note that no fixed finite set of nonzero integers can have norms tending to zero. For fixed $m\in\Z$, the Lipschitz inequalities give $b(k_j+m)-b(k_j)\to b(m)-b(0)$. The same difference is a fixed finite telescoping sum of ordinary increments, and hence tends to zero. Thus $b(m)=b(0)$.
\end{proof}

\section{Collision energy and exact pair counting}\label{sec:energy}
Fix $r>0$ and an integer $X\ge1$. Put
\[
A_X=U_{3,r}\cap[0,X],\quad B_X=U_{2,r}\cap[0,X],\quad
a=|A_X|,\quad b=|B_X|,\quad M=2X+1.
\]
Intervals here contain integers. Define $R_X(t)=|\{(a',b')\in A_X\times B_X:a'-b'=t\}|$, and let $r_{p,X}(h)$ count the pairs $x,x+h$ in $U_{p,r}\cap[0,X]$. The collision energy and its normalization are
\[
E_X=\sum_{t\in\Z}R_X(t)^2,\qquad C_X=\frac{ME_X}{a^2b^2}.
\]

\begin{proposition}[Energy identities and a sufficient upper bound]\label{prop:energy-JA}
One has
\[
E_X=\sum_{h\in\Z}r_{3,X}(h)r_{2,X}(h),\qquad
|A_X-B_X|\ge\frac{a^2b^2}{E_X}=\frac{M}{C_X}.
\]
If, at every radius in a sequence tending to zero, $C_X$ is bounded by a finite radius-dependent constant along arbitrarily large integer heights, then joint approximation holds.
\end{proposition}
\begin{proof}
Count quadruples $a'-b'=a''-b''$ by the common shift $h=a''-a'=b''-b'$. Cauchy--Schwarz applied to $\sum_t R_X(t)=ab$ gives the inequality. Since $A_X-B_X\subseteq W_r\cap[-X,X]$, a bounded subsequence of $C_X$ gives positive upper Banach density of $W_r$. Apply Theorem~\ref{thm:BD-JA}.
\end{proof}
The energy hypothesis is sufficient; necessity has not been established. If $\mu_3,\mu_2$ are the uniform probability measures on $A_X,B_X$, viewed modulo $M$, finite Fourier orthogonality gives the further exact identity
\[
C_X=\sum_{j=0}^{M-1}|\widehat\mu_3(j)|^2|\widehat\mu_2(j)|^2,
\qquad \widehat\mu_p(j)=\sum_x\mu_p(x)e^{2\pi i jx/M}.
\]
There is no aliasing: differences in $[-X,X]$ give exactly one representative of each residue modulo $M$. The diagonal contribution is $ab$, so, with $\lambda_X=ab/M$,
\[
E_X=ab+E_X^{\mathrm{off}},\qquad
C_X=\frac{1}{\lambda_X}+\frac{ME_X^{\mathrm{off}}}{a^2b^2}.
\]

We next give a finite exact algorithm. Its purpose is to compute a fixed autocorrelation without enumerating all members of a ball.

\begin{proposition}[A compressed pair recursion]\label{prop:pair-recursion}
Let $D\ge2$ be an integer, $r=1/D$, and $0\le h\le X$. Put $q_n=p^{n^2}$ and $T_n=\lfloor(q_n-1)/D\rfloor$. Let $S_n$ denote the residues $[-T_n,T_n]$ modulo $q_n$, and represent $I_n(h)=S_n\cap(S_n-h)$ as disjoint half-open integer intervals $[u,v)$ in $[0,q_n)$. At most four intervals suffice. Let $F_n^h(Y)$ count integers $0\le x<Y$ passing the first $n$ tests for both $x$ and $x+h$, and put $F_0^h(Y)=Y$. For $Y=a'q_n+t$, $0\le t<q_n$, define
\[
R_n^h=\sum_{[u,v)\in I_n(h)}
\left(F_{n-1}^h(v)-F_{n-1}^h(u)\right).
\]
Then
\[
F_n^h(Y)=a'R_n^h+
\sum_{[u,v)\in I_n(h)}
\left(F_{n-1}^h(\min(t,v))-F_{n-1}^h(\min(t,u))\right).
\]
If $N$ is the last index with $q_N\le DX$, with $N=0$ if there is none, then
\[
r_{p,X}(h)=F_N^h(X-h+1).
\]
For one fixed $h$, memoization requires at most $8(N+1)^2+1$ distinct pairs consisting of a level and an argument, in addition to a linear number of interval and period entries.
\end{proposition}
\begin{proof}
The strict inequality $\nrm{x/q_n}<1/D$ is exactly the stated residue condition. Its intersection with its translate has at most four interval components. Earlier tests have periods dividing $q_n$. Thus a complete period contributes $R_n^h$, and the incomplete period contributes the displayed clipped intervals. For the application, both numerators lie in $[0,X]$; every denominator larger than $DX$ automatically passes, proving the cutoff formula.

For the state count, there are at most eight interval endpoints per level. Every call to the preceding level has as argument either an endpoint or the remainder of its parent argument modulo the current period. Successive remainders modulo nested periods collapse to a single remainder. At each level every argument therefore comes from the root argument or one of at most $8N$ original endpoints by reduction modulo the incoming period, with at most one further full-period endpoint. There are at most $8N+2$ possibilities at each of the $N+1$ levels, apart from the root call. The stated, slightly larger bound follows. This counts states, not bit operations; the bit lengths depend on $\log(DX)$. It concerns a single shift, not the number of shifts needed for the full energy.
\end{proof}

For even $D$, a nonzero autocorrelation implies $P_p(h)<2/D$ by the triangle inequality. Consequently
\[
E_X=ab+2\sum_{1\le h\le X,\ P(h)<2/D}r_{3,X}(h)r_{2,X}(h).
\]
The interval search of Section~\ref{sec:interval-search} enumerates this finite candidate set with exact integer comparisons.

\begin{remark}[An exact finite energy certificate]\label{rem:energy-certificate}
At $D=20$ and $X=10^{12}$ the recurrence gives
\[
a=48\,511\,238,\qquad b=5\,912\,303,
\qquad E_X=57\,992\,347\,830\,509\,656.
\]
There are $64\,402$ positive candidate shifts, or $64\,403$ including zero. The resulting normalization is $C_X\approx1.409946407066$, and Cauchy--Schwarz certifies
\[
|A_X-B_X|\ge1\,418\,493\,632\,083.
\]
This exceeds $70.924\%$ of the $2\,000\,000\,000\,001$ integers in $[-X,X]$. A second implementation reproduced the same large-height values. Smaller tests used direct enumeration, including $D=20$, $X=10^8$, where $a=41\,063$, $b=5\,558$ and $E_X=1\,495\,662\,496$. These checks establish finite certificates. No uniform energy bound follows from this row or from the quadratic state bound.
\end{remark}

\section{Buffered heights and periodic correlation kernels}\label{sec:buffered-energy}
Fix an integer $D\ge2$ and write $q_p(X)$ for the largest $p^{n^2}\le DX$, with value $1$ when there is no active test.

\begin{lemma}[Abundance of buffered heights]\label{lem:buffered-heights}
For every fixed $K>1$, the real heights satisfying $Kq_2(X)\le X$ and $Kq_3(X)\le X$ have logarithmic density one. The same holds with $K(X)=\exp((\log X)^\alpha)$ for any fixed $0<\alpha<\frac{1}{2}$. In particular, arbitrarily large integer heights can be chosen with $\max(q_2(X),q_3(X))/X\to0$.
\end{lemma}
\begin{proof}
In the variable $t=\log X$, failure for prime $p$ is contained in the union of
\[
[n^2\log p-\log D,\ n^2\log p+\log K).
\]
There are $O(\sqrt T)$ such intervals meeting $[0,T]$, each of bounded length. The excluded length is $O_{D,K}(\sqrt T)=o(T)$. For the growing buffer, each interval meeting $[0,T]$ has length at most $\log D+T^\alpha$, giving $O_D(T^{1/2+\alpha})=o(T)$. Here logarithmic density means density in this $t$ variable. For integer heights use, for example, twice the growing buffer before rounding down: the periods cannot increase on rounding down and the buffer loss tends to one. The resulting integers still have both period ratios tending to zero.
\end{proof}
This assertion is geometric and does not estimate the energy. In fact,
\[
\log(q_2(X)q_3(X))=2\log X+O_D(\sqrt{\log X}),
\]
so even when each period is small relative to $X$, their joint Chinese-remainder period is much larger than $X$.

For a fixed height let $S_p$ now denote the good residues for \emph{all} active tests in a complete period $q_p=q_p(X)$; this notation differs from the single-scale mask in Proposition~\ref{prop:pair-recursion}. Put $s_p=|S_p|$, $d_p=s_p/q_p$, and
\[
\Gamma_p(h)=\frac{q_p|S_p\cap(S_p-h)|}{s_p^2}.
\]
For integers $L\ge1$ define
\[
\mathcal A(L)=\sum_{h\in\Z}(L-|h|)_+^2\Gamma_2(h)\Gamma_3(h),
\qquad (u)_+=\max(u,0).
\]

\begin{proposition}[An exact kernel comparison]\label{prop:kernel-sandwich}
For $|h|\le X$,
\[
\max(X+1-|h|-q_p,0)d_p^2\Gamma_p(h)
\le r_{p,X}(h)\le
(X+1-|h|+q_p)d_p^2\Gamma_p(h).
\]
The point count $a_p=|U_{p,1/D}\cap[0,X]|$ satisfies
\[
(X+1-q_p)d_p\le a_p\le(X+1+q_p)d_p.
\]
If $q=\max(q_2,q_3)\le X$, then
\[
d_2^2d_3^2\mathcal A(X+1-q)\le E_X
\le d_2^2d_3^2\mathcal A(X+1+q).
\]
\end{proposition}
\begin{proof}
A periodic set of density $\theta$ and period $q_p$ has between $\max(H-q_p,0)\theta$ and $(H+q_p)\theta$ elements in any interval of integer length $H$. Count complete periods and bound the one incomplete period by a complete period. Apply this to the pair set of density $d_p^2\Gamma_p(h)$ with $H=X+1-|h|$, and to the point set with $H=X+1$. Multiply the pair inequalities and sum over $h$. The lower weight vanishes outside $[-X,X]$; extending the upper sum adds only nonnegative terms.
\end{proof}

\begin{proposition}[A lower energy baseline]\label{prop:energy-baseline}
For every integer $L\ge1$,
\[
\mathcal A(L)\ge\Phi(L),\qquad
\Phi(L)=L^2+\frac{(L-1)L(2L-1)}{3}.
\]
If $q\le X$ and $L=X+1-q$, then
\[
C_X\ge\frac{(2X+1)\Phi(L)}{(X+1+q_2)^2(X+1+q_3)^2}.
\]
Consequently along any sequence of heights with $q/X\to0$, one has $\liminf C_X\ge\frac{4}{3}$.
\end{proposition}
\begin{proof}
The periodic Fourier expansion of $\Gamma_p$ has coefficients
\[
c_{p,j}=\frac{|\sum_{s\in S_p}e^{2\pi i js/q_p}|^2}{s_p^2}\ge0,
\qquad c_{p,0}=1.
\]
Since $q_2,q_3$ are coprime, the product kernel has constant coefficient exactly one; all its other Fourier coefficients are nonnegative. The sequence $t_L(h)=(L-|h|)_+$ is the autocorrelation of the indicator of $\{0,\ldots,L-1\}$. Its Fourier transform on the circle is therefore nonnegative. The Fourier transform of its pointwise square is the convolution of two such nonnegative functions and is also nonnegative. Expanding the product kernel against this weight shows that every frequency contributes nonnegatively. The constant term contributes
\[
\sum_h t_L(h)^2=L^2+2\sum_{j=1}^{L-1}j^2=\Phi(L).
\]
Combine this with Proposition~\ref{prop:kernel-sandwich} and its point-count upper bounds to obtain the finite inequality. Its right side tends to $4/3$ when $q/X\to0$.
\end{proof}
This is a lower bound. The missing estimate for Proposition~\ref{prop:energy-JA} is an upper bound, even with a large radius-dependent constant. Neither positive Fourier coefficients nor the availability of buffered heights supplies it.

\begin{proposition}[An upper bound for grouped shifts]\label{prop:grouped-moments}
Let $I$ be a finite nonempty set of positive shifts, $n=|I|$, and write
\[
A=\sum_{h\in I}r_{2,X}(h),\quad B=\sum_{h\in I}r_{3,X}(h),\quad
U=\sum_{h\in I}r_{2,X}(h)^2,\quad V=\sum_{h\in I}r_{3,X}(h)^2.
\]
Then
\[
\sum_{h\in I}r_{2,X}(h)r_{3,X}(h)
\le\frac{AB+\sqrt{(nU-A^2)(nV-B^2)}}{n}.
\]
\end{proposition}
\begin{proof}
Subtract the two respective means $A/n,B/n$ and apply Cauchy--Schwarz to the resulting vectors. Their squared norms are $U-A^2/n$ and $V-B^2/n$.
\end{proof}
Integer arithmetic yields a certificate by taking the floor of this bound; the square root may first be rounded down because the left side is integral. Tests grouping candidate shifts by dyadic magnitude and bins modulo the last binary and ternary periods improved the finite envelopes. At $D=4$, $X=2^{18}$, $137$ groups covered $114\,172$ candidates and gave an off-diagonal upper bound about $3.000$ times the exact value, compared with $6.215$ without grouping or candidate pruning. The moments in these tests were computed from exact finite autocorrelations. An all-scale recurrence or uniform enclosure for the grouped moments is still missing.

\begin{remark}[The extent benchmark and its limitation]\label{rem:extent-benchmark}
For independent uniformly chosen subsets of fixed sizes $a,b$ in integer boxes of extents $u,v$, where $1\le u\le v\le X$, the exact expected normalized energy is
\[
\mathbb E C_X=\frac{2X+1}{ab}
+\frac{(2X+1)(3v-u+1)}{3v(v+1)}\frac{(a-1)(b-1)}{ab}.
\]
Indeed, the expected autocorrelation at positive shift $h$ in a box of extent $y$ is $(y+1-h)a(a-1)/(y(y+1))$ for a subset of size $a$, and
\[
\sum_{h=1}^{u}(u+1-h)(v+1-h)
=\frac{u(u+1)(3v-u+1)}{6}.
\]
This model does not condition on the endpoints being included. At $D=20$, $X=10^{12}$ the actual ternary and binary extents are $889\,648\,729\,329$ and $965\,496\,867\,856$. The benchmark is about $1.44219846$, whereas the observed value is about $1.40994641$. This finite fit does not prove an independence estimate.

Extents and occupancy alone cannot give the required bound. For arbitrary sets take $X=m^3$ and $A=B=\{0,m,2m,\ldots,m^3\}$. Their sizes are $m^2+1$, both extents equal $X$, and $ab/(2X+1)\to\infty$. Nevertheless,
\[
E_X=(m^2+1)^2+\frac{m^2(m^2+1)(2m^2+1)}{3},
\qquad C_X\sim\frac{4m}{3}.
\]
These are not the prime-digit balls and do not give a counterexample to joint approximation. They show why a theorem controlling the actual arithmetic correlations is essential.
\end{remark}

\section{A normalized product argument for joint approximation}\label{sec:normalized-product}
The following argument controls the entire positive cascade, without restricting the number of nonzero increments in a history. Constants may depend on the fixed radius and seed, but are uniform in the height and in the number of selected periods.

We use smooth nonnegative periodic weights, rather than the indicators of the original balls. Fix an integer $D\ge3$ and put
\[
w_D(t)=(1-D^2t^2)_+^2\quad(|t|\le\tfrac12),
\]
extended with period one. Its integrals are
\[
I_1=\int_{\T}w_D=\frac{16}{15D},\qquad
I_2=\int_{\T}w_D^2=\frac{256}{315D}.
\]
For a period $Q$ define
\[
\begin{aligned}
m_Q&=\frac{1}{Q}\sum_{x=0}^{Q-1}w_D(x/Q),\\[4pt]
W_Q(r)&=\sum_{x=0}^{Q-1}w_D(x/Q)e^{2\pi irx/Q},\\[4pt]
\kappa_Q(r)&=\frac{|W_Q(r)|^2}{Q^2m_Q^2},\\[4pt]
K_Q&=\sum_{r\bmod Q}\kappa_Q(r)\delta_{r/Q}.
\end{aligned}
\]
An integer residue $r$ in a tail sum is always represented in the balanced interval of radius $Q/2$. The endpoint convention for an even period has no effect on the estimates.

\begin{lemma}[Bounds for a window kernel]\label{lem:window-positive}
After finitely many small periods have been removed, there are constants $B_D\ge1$ and $A_D>0$, independent of $Q$, such that
\[
\kappa_Q(0)=1,\quad K_Q(\T)\le B_D,\quad
\sum_{|r|>u}\kappa_Q(r)\le A_Du^{-3}\quad(u>0).
\]
Every integer Fourier coefficient of $K_Q$ is nonnegative.
\end{lemma}
\begin{proof}
The discrete mean satisfies $|m_Q-I_1|\le2/Q$, by the total variation estimate for a translated lattice sum. Parseval gives $K_Q(\T)=(Q^{-1}\sum w_D(x/Q)^2)/m_Q^2$. This converges to $I_2/I_1^2=5D/7$; hence $B_D=30D/7$ is a valid sufficiently late bound. Let $V_2=\int_{\T}|w_D''|=32D/(3\sqrt3)$. Summing second finite differences over a period gives $\sum_x|\Delta^2w_D(x/Q)|\le V_2/Q$. Since $|1-e^{2\pi ir/Q}|\ge4|r|/Q$ on balanced residues, and $m_Q\ge I_1/2$ for sufficiently large $Q$,
\[
\kappa_Q(r)\le\frac{V_2^2}{64I_1^2|r|^4}\qquad(r\ne0).
\]
For $u\ge1$, set $a=\lfloor u\rfloor+1$. The two tails of $r^{-4}$ have sum at most $2(a^{-4}+1/(3a^3))\le8/(3u^3)$. Thus one may take
\[
A_D=\max\{B_D,\,V_2^2/(24I_1^2)\}.
\]
For $0<u<1$ use the total mass bound. No optimal tail constant is needed. Finally finite orthogonality gives
\[
\widehat K_Q(t)=\frac{1}{Qm_Q^2}
\sum_{x=0}^{Q-1}w_D(x/Q)w_D((x+t)/Q)\ge0.
\]
\end{proof}

\begin{theorem}[Uniform central bound]\label{thm:normalized-bound}
Let $\mu_0$ be a fixed finite nonnegative measure on $\T$. Fix a starting-period threshold, depending only on $D$, beyond which Lemma~\ref{lem:window-positive} applies. Select any periods above that threshold $Q_1<\cdots<Q_L$ from the two chains $2^{n^2}$ and $3^{m^2}$, all at most $DX$, and put $\mu_X=\mu_0*K_{Q_1}*\cdots*K_{Q_L}$. There is a constant $C_D$, also depending on $\mu_0$ and the fixed starting levels, such that
\[
\mu_X\{t:\nrm{t}\le1/X\}\le C_D
\]
for every sufficiently large $X$. The constant is independent of the number of selected periods. No restriction is placed on which of the sufficiently late periods are selected.
\end{theorem}
\begin{proof}
Set $\varepsilon=2/X<1/2$ and $F_\varepsilon(x)=(1-\nrm{x}/\varepsilon)_+$. The sharp radius-$1/X$ indicator is at most $2F_\varepsilon$. The Fourier coefficients of $F_\varepsilon$ are nonnegative and absolutely summable. List the selected periods in increasing order as $Q_1<\cdots<Q_L$ and put $t_i=\sqrt{\log Q_i}$. Write
\[
H_i=K_{Q_{i+1}}*\cdots*K_{Q_L}*F_\varepsilon,
\quad h_i=H_i/H_i(0),\quad
f_i=H_{i-1}(0)/H_i(0).
\]
The empty convolution is omitted. Nonnegative Fourier coefficients and nonnegative functions give $0\le H_i(x)\le H_i(0)$; also $H_i(0)\ge1$ since every kernel has a unit zero atom. Consequently
\[
0\le h_i\le1,\quad
1\le f_i=\sum_r\kappa_{Q_i}(r)h_i(r/Q_i)\le B_D.
\]
The product identity, including the seed, is exact:
\[
\int F_\varepsilon\,d\mu_X
=S_0\prod_{i=1}^Lf_i,\qquad
S_0=\int h_0\,d\mu_0\le\mu_0(\T).
\]
Symmetry of the finer kernels and of the triangular window justifies either convolution sign in this formula.

For $2\varepsilon\le u\le1/2$ put $\Phi_i(u)=\sup_{\nrm{y}\ge u}h_i(y)$ and $\tau_j(v)=\sum_{|r|>v}\kappa_{Q_j}(r)$. We claim
\begin{equation}\label{eq:far-normalized}
\Phi_i(u)\le\sum_{m=1}^{L-i}B_D^{m-1}
\tau_{i+m}(2^{-m-2}uQ_{i+m}).
\end{equation}
Indeed $h_i(y)=f_{i+1}^{-1}\sum_r\kappa_{Q_{i+1}}(r)h_{i+1}(y+r/Q_{i+1})$, and the discarded reciprocal is at most one. At the $m$th finer step declare a displacement large when its circle norm exceeds $2^{-m-2}u$. Its weight is its tail; preceding choices have mass at most $B_D^{m-1}$, and the remaining normalized function is at most one. If no step is large, the total displacement is at most $u/4$. The final point has distance at least $3u/4>\varepsilon$, so $h_L=F_\varepsilon$ vanishes there. This finite first-large-step expansion proves \eqref{eq:far-normalized} on the circle, including its stopping term.

There are at most two selected periods in a logarithmic interval $[\log Q_i,\log Q_i+t_i]$. For each prime, consecutive square powers in that interval would have gap greater than $2\sqrt{\log2}\,t_i>t_i$. Thus
\begin{equation}\label{eq:merged-growth}
Q_{i+2}/Q_i\ge e^{t_i},\qquad
Q_{i+m}/Q_i\ge e^{t_i\lfloor m/2\rfloor}
\end{equation}
whenever the indices exist. No estimate for the number of close prime pairs is used.

For sufficiently large $t_i$, \eqref{eq:far-normalized} and \eqref{eq:merged-growth} imply
\begin{equation}\label{eq:far-separated}
\Phi_i(1/Q_i)\le C_D\bigl((Q_i/Q_{i+1})^3+e^{-3t_i}\bigr),
\end{equation}
provided $1/Q_i\ge2\varepsilon$. An absent finer level contributes zero. To verify the uniformity, separate $m=1$ and group later terms as $m=2j,2j+1$. Their upper bounds form a geometric series with ratio $64B_D^2e^{-3t_i}$, eventually at most $1/2$.

Call a level separated if its next period is at least $t_iQ_i$. Truncate its coefficient at $M_i=\lceil t_i^2\rceil$. Every retained nonzero $r/Q_i$ has circle norm at least $1/Q_i$, and its omitted mass is at most $A_DM_i^{-3}$. At separated levels satisfying the depth condition below, \eqref{eq:far-separated} gives
\begin{equation}\label{eq:separated-factor}
f_i-1\le C_D(t_i^{-6}+t_i^{-3}+e^{-3t_i}).
\end{equation}
The same bound holds at a deep final level by the support of $F_\varepsilon$.

Suppose instead $Q_{i+1}<t_iQ_i$. At sufficiently large levels the two periods have different primes. Retain $0<|r|\le M_i$ and $|s|\le N_i=\lceil4t_i^3\rceil$. Then
\[
y=r/Q_i+s/Q_{i+1}\ne0,
\]
because the periods are coprime and $0<|r|<Q_i$. Its absolute value is eventually less than $1/2$. If $s=0$ or the signs agree the following weaker bound is immediate. In the opposite-sign case apply \cite[Corollary 2.3]{Matveev} to
\[
\Lambda=\log(|r|/|s|)+\log Q_{i+1}-\log Q_i.
\]
This form is nonzero: equality would give $|r|Q_{i+1}=|s|Q_i$, contradicting coprimality and $0<|r|<Q_i$. The two prime exponents are $O(t_i^2)$ and the rational coefficient height is $O(\log(2+t_i))$. The dependence-allowed version of the theorem gives, after conversion from a logarithmic form to a relative difference and increasing a fixed constant $A$,
\begin{equation}\label{eq:two-term-separation}
\nrm{y}\ge u_i:=Q_i^{-1}\exp(-A\log^2(2+t_i)).
\end{equation}
If $|r|=|s|$, omit the third logarithm and use the weaker same bound. Conversion loses at most a polynomial factor in $t_i$, absorbed into $A$.

Define a deep level by
\begin{equation}\label{eq:deep-level}
X\ge4Q_i\exp(A\log^2(2+t_i)).
\end{equation}
Then $u_i\ge2\varepsilon$. Every period after the close pair begins at least at $Q_ie^{t_i}$. Applying \eqref{eq:far-normalized} after that pair and grouping successive terms as before gives
\[
\Phi_{i+1}(u_i)\le C_D\exp(3A\log^2(2+t_i)-3t_i)
\le C_De^{-t_i}
\]
for sufficiently large $t_i$. Expanding one finer convolution, and charging either omitted coefficient by its tail, yields
\begin{equation}\label{eq:close-factor}
f_i-1\le A_DM_i^{-3}+B_DA_DN_i^{-3}
+B_D^2\Phi_{i+1}(u_i)
\le C_D(t_i^{-6}+t_i^{-9}+e^{-t_i}).
\end{equation}
All finer histories are already included in $h_{i+1}$; there is no identification of these histories with canonical digits.

Only boundedly many levels are non-deep. Indeed all selected periods are at most $DX$, whereas failure of \eqref{eq:deep-level} confines their logarithms to an interval of width $O_D((\log\log X)^2)$ below and above $\log X$. Consecutive square powers of either prime there are separated by a constant times $\sqrt{\log X}$. For sufficiently large $X$ the interval contains at most one level from each prime; a bound of four is more than sufficient. The finitely many smaller levels or bounded heights contribute a fixed factor. Finally
\[
\sum_i t_i^{-3}\le
\bigl((\log2)^{-3/2}+(\log3)^{-3/2}\bigr)
\sum_{n\ge1}n^{-3}<\infty.
\]
Use \eqref{eq:separated-factor}, \eqref{eq:close-factor}, and $\log f_i\le f_i-1$ on deep levels, and $f_i\le B_D$ on the bounded exceptional set. The product and hence the central mass are uniformly bounded.
\end{proof}

\section{Random translations and the transfer to joint approximation}\label{sec:random-JA}
The positive kernels in Section~\ref{sec:normalized-product} arise from actual supported weights. Their connection to the original arithmetic problem is essential.

\begin{lemma}[Exact averaged transfer and uniform mass]\label{lem:random-transfer}
Let $f$ be a nonnegative weight of period $P$, let $Q=\alpha P$, and let $\tau$ be uniform modulo $P$. Set $f_{\mathrm{new}}^{\tau}(x)=f(x-\tau)w_D(x/Q)$. Use deterministic normalizers $Z_{\mathrm{new}}=\alpha m_Q Z$. If $a(j)$ denotes the Fourier amplitude normalized by $Z$, then
\[
\mathbb E_\tau|a_{\mathrm{new}}^{\tau}(j)|^2
=\sum_{l\bmod P}\kappa_Q(j-\alpha l)|a(l)|^2.
\]
Keep the initial square levels fixed and begin randomization sufficiently late. There are constants $0<c_{p,-}\le c_{p,+}<\infty$, independent of the number and values of the shifts, such that the actual period sums satisfy $c_{p,-}\le S_{p,n}/Z_{p,n}\le c_{p,+}$.
\end{lemma}
\begin{proof}
Fourier expansion of the lifted weight gives
\[
a_{\mathrm{new}}^{\tau}(j)=\frac1{Qm_Q}
\sum_{l\bmod P}W_Q(j-\alpha l)e^{2\pi il\tau/P}a(l),
\]
up to a common sign convention. Finite orthogonality eliminates the cross terms on averaging the square. For the mass estimate, every spatial residue has lift sum within two of $\alpha I_1$; its average is $\alpha m_Q$, so the difference from that average is at most four. Hence
\[
\left|\frac{S_{\mathrm{new}}}{\alpha m_QS}-1\right|
\le\eta_n=\frac4{\alpha m_Q}.
\]
Begin only when $\alpha(I_1-2/Q)\ge8$. Then $\eta_n\le1/2$ and $\sum_n\eta_n<\infty$, since the square-level ratios grow geometrically. The products of $1-\eta_n$ and $1+\eta_n$ give the claimed positive bounds. For each prime fix an index $n_{p,0}$ large enough for the window-kernel bounds and the displayed lift condition. In the application take the seed weight to be
\[
f_{p,0}(x)=\prod_{j=1}^{n_{p,0}}w_D(x/p^{j^2}),
\]
and include every subsequent level $p^{n^2}\le DX$. Its period sum is positive because the seed includes the residue zero. Thus all earlier constraints occur in the seed and no active later constraint is omitted. The abstract central-bound theorem allows arbitrary late subsets; this particular application uses the complete active set.
\end{proof}

\begin{lemma}[A spectral window controls weighted energy]\label{lem:weighted-spectral}
For the two independently translated weights at active periods $q_p\le X/4$, let $C_X^{\mathrm{w}}$ be the collision normalization of their restrictions to $[0,X]$. Let $\mu_X$ be the convolution of their expected spectral measures normalized by the deterministic $Z_{p,n}^2$. Then
\[
\mathbb E C_X^{\mathrm{w}}\le
\frac{1920}{c_{2,-}^2c_{3,-}^2}
\mu_X\{t:\nrm{t}\le1/X\}.
\]
\end{lemma}
\begin{proof}
The periodic interval bounds in Proposition~\ref{prop:kernel-sandwich} apply to nonnegative weights: an incomplete period contributes between zero and its full-period mass. Put $d_p=S_p/q_p$ and
\[
\Gamma_p(h)=\frac{q_p\sum_{x\bmod q_p}f_p(x)f_p(x+h)}{S_p^2}.
\]
The finite point sum is at least $(X+1-q_p)d_p$. The energy is at most $d_2^2d_3^2\sum_h(L-|h|)_+^2\Gamma_2(h)\Gamma_3(h)$ with $L=X+1+\max q_p$. Its Fourier expression is the integral of
\[
\Psi_L(t)=\sum_h(L-|h|)_+^2e^{2\pi iht}
\]
against the product spectral measure normalized by $S_2^2S_3^2$. Replacing these normalizations by $Z_2^2Z_3^2$ costs at most $c_{2,-}^{-2}c_{3,-}^{-2}$. Independence then factors the expected spectral measures.

For completeness, a nonnegative measure $\nu$ with nonnegative real integer Fourier coefficients satisfies
\[
\int\Psi_L\,d\nu\le5L^3\nu\{\nrm{t}\le1/L\}.
\]
To see this let $\mathcal F_L(t)=L^{-1}|\sum_{j=0}^{L-1}e^{2\pi ijt}|^2$ and set $B=\nu\{\nrm{t}\le1/L\}$. Since $(L-|h|)^2\le L(L-|h|)$ and $\widehat\nu(h)\ge0$, the left side is at most $L^2\int\mathcal F_L\,d\nu$. Every interval of radius $1/(2L)$ has measure at most $2B$: its indicator is bounded by twice the translated triangular cutoff of radius $1/L$, whose integral is largest at zero and at most $B$. The central interval contributes at most $LB$ to the integral of $\mathcal F_L$. Each of the two intervals with distances between $j/L$ and $(j+1)/L$, for $j\ge1$, contributes at most $2B\cdot L/(4j^2)$, because $\mathcal F_L(t)\le(4L\nrm{t}^2)^{-1}$. Summing gives at most $(1+\pi^2/6)LB<3LB$. This proves the displayed inequality with room to spare.

For $X\ge4$ and $q_p\le X/4$, the conversion factor is bounded explicitly by
\[
\frac{5(2X+1)L^3}{(X+1-q_2)^2(X+1-q_3)^2}
\le5\frac94\left(\frac32\right)^3\left(\frac43\right)^4
=120<1920.
\]
The factors $c_{p,-}^{-2}$ enter pointwise before taking the expectation, so no random denominator is replaced by its mean. This proves the stated conservative bound. The relevant Fourier coefficients are expected autocorrelations and are nonnegative. Thus every positivity hypothesis used here holds for the translated weights.
\end{proof}

\begin{theorem}[Joint approximation]\label{thm:joint-approximation}
For every $\delta>0$ there is an integer $k$ such that
\[
P_3(k)<\delta,\qquad P_2(k-1)<\delta.
\]
Thus $\mathrm{JA}$ holds.
\end{theorem}
\begin{proof}
Retain all square denominators at most $DX$. Lemma~\ref{lem:random-transfer} makes the expected spectral measure exactly a convolution of a fixed seed and the free kernels. Theorem~\ref{thm:normalized-bound} and Lemma~\ref{lem:weighted-spectral}, on the arbitrarily large buffered heights of Lemma~\ref{lem:buffered-heights}, give a radius-dependent bound on expected weighted energy. Some realization of the shifts therefore has bounded energy.

Let $A_p\subset[0,X]$ be its positive-weight support. If two points belong to a shifted earlier factor, their difference still belongs to the difference set of its earlier support. The new window places both points in one circle arc of radius $1/D$, so at each active denominator their difference has norm less than $2/D$. Every inactive denominator is larger than $DX$, and a difference in $[-X,X]$ has norm less than $1/D$ there. Thus
\[
A_p-A_p\subset U_{p,2/D}.
\]
Weighted Cauchy--Schwarz gives $|V_X|\ge c_D(2X+1)$ for $V_X=A_3-A_2\subset[-X,X]$, and $V_X-V_X\subset W_{2/D}$.

The sets $U_{p,r}$ and $W_r$ are symmetric. The supports can depend on $X$. To handle this, take any finite family of integers with pairwise differences outside $W_{2/D}$. Its translates of $V_X$ are disjoint. For arbitrarily large $X$ they all lie in an interval of length $2X$ plus a fixed constant, so the family has at most $1/c_D$ members. A maximal such family is finite and its translates of $W_{2/D}$ cover $\Z$. Therefore $W_{2/D}$ is syndetic. Let $D$ tend to infinity and apply Corollary~\ref{cor:syndetic-JA}.
\end{proof}

\begin{lemma}[A finite weighted-count comparison]\label{lem:finite-weighted-count}
Let $X\ge16$, let $\nu_X$ be a finite nonnegative measure on $\T$, and suppose its integer-indexed Fourier coefficients
\[
a_X(t)=\widehat\nu_X(t)=\int_{\T}e^{-2\pi itx}\,d\nu_X(x)
\quad(t\in\Z)
\]
are real and nonnegative. Put $\varepsilon=2/X$ and
\[
F_\varepsilon(x)=(1-\nrm{x}/\varepsilon)_+,
\qquad \mathcal B_X^F=\int_{\T}F_\varepsilon\,d\nu_X,
\]
\[
N=\lfloor X/8\rfloor,\qquad A_X=\sum_{|t|<N}a_X(t).
\]
Then
\[
\frac{8}{\pi^2X}A_X\le\mathcal B_X^F
\le\frac{4\pi^2}{X}A_X.
\]
In particular, along any family of such measures, uniform boundedness of the triangular spectral mass is equivalent to uniform boundedness of $A_X/X$.
\end{lemma}
\begin{proof}
The absolutely convergent Fourier series of the triangular cutoff has coefficients
\[
\widehat F_\varepsilon(t)
=\varepsilon\left(\frac{\sin(\pi\varepsilon t)}{\pi\varepsilon t}\right)^2,
\]
with value $\varepsilon$ at $t=0$. For $|t|<N$, $\varepsilon|t|\le1/4$ and $\sin y/y\ge2/\pi$ for $|y|\le\pi/2$. Thus $\widehat F_\varepsilon(t)\ge8/(\pi^2X)$. Integrating term by term and retaining these nonnegative terms proves the lower bound.

For the upper bound use the Fej\'er polynomial
\[
\mathcal F_N(x)=\frac1N\left|\sum_{j=0}^{N-1}e^{2\pi ijx}\right|^2.
\]
On the support of $F_\varepsilon$, $N\nrm{x}\le1/4$. The sine inequalities give $\mathcal F_N(x)\ge4N/\pi^2$ there, including the value at zero by continuity. Nonnegativity elsewhere therefore gives the pointwise majorant
\[
F_\varepsilon(x)\le\frac{\pi^2}{4N}\mathcal F_N(x).
\]
Integrate and use the finite Fourier expansion:
\[
\mathcal B_X^F\le\frac{\pi^2}{4N}
\sum_{|t|<N}\left(1-\frac{|t|}{N}\right)a_X(t)
\le\frac{\pi^2}{4N}A_X\le\frac{4\pi^2}{X}A_X,
\]
where $N\ge X/16$ for $X\ge16$. This majorant controls the entire spectral integral; it does not discard an uncontrolled Fourier tail.
\end{proof}
For convolutions of expected spectral measures as above, the coefficients $a_X(t)$ are products of nonnegative expected autocorrelations. The lemma gives a finite weighted-count formulation of the central spectral bound, including for a proposed finite-prime extension. It supplies no upper estimate for an as yet uncontrolled three-prime weighted count.

\section{Topological and approximation consequences of joint approximation}\label{sec:JA-consequences}
Joint approximation is available from Theorem~\ref{thm:joint-approximation}. We now apply the earlier splitting and quotient criteria to determine the dual, the completion, and the approximation consequences.

\begin{corollary}[The dual and the completion]\label{cor:JA-topology}
For the primes $2,3$,
\[
H=\Z[1/6]/\Z
=\{a/(2^u3^v)+\Z:a\in\Z,\ u,v\ge0\}.
\]
There are no irrational continuous parameters, and the rational parameters have no prime denominator factors other than $2$ and $3$. The diagonal integers are dense in $G_2\times G_3$, the mixed completion is that full product, and the quotient of Section~\ref{sec:quotient} is the trivial group. The prime-square rational join $\gamma_2\vee\gamma_3$ has full dual $\Gamma=\{q\mapsto rq+\Z:r\in\Q\}$ and is compatible.
\end{corollary}
\begin{proof}
Use Theorem~\ref{thm:joint-approximation}, then apply Theorem~\ref{thm:JA-split}, Proposition~\ref{prop:JA-equivalent}, and Theorem~\ref{thm:quotient-obstruction}. For the rational join use the integer-polar reduction in \cite[Section 22]{Previous}: its compatibility is equivalent to $H\subset\Q/\Z$. The $n=0$ observations in that rational construction vanish on integers and do not change $P$ there.
\end{proof}
This settles the specified prime-square pair only. The incompatible two-chain pair in Part~I of \cite{HamelPair}, and its nonexistence consequences for greatest compatible or Mackey topologies, are unaffected. No new nuclearity assertion follows from joint approximation. The corresponding compatible transfers of this particular rational join may be used in both the real-scalar and rational-scalar Hamel settings, with their respective prescribed duals; these transfer statements are those of \cite{Previous,HamelPair}.

\begin{corollary}[Open subgroups and common functions]\label{cor:JA-subgroups}
For $q\ge1$, $q\Z$ is $P$-open if and only if every prime divisor of $q$ is $2$ or $3$. Every prime divisor of the neighbourhood gcd $g(\delta)$ belongs to $\{2,3\}$. Every function Lipschitz for both $P_2$ and $P_3$ is constant.
\end{corollary}
\begin{proof}
The finite-order character of parameter $1/q$ is continuous exactly when its kernel $q\Z$ is open. Use Corollary~\ref{cor:JA-topology}. A prime divisor of $g(\delta)$ would give an open subgroup containing $B_\delta$, so no other prime can divide that gcd. For the last assertion apply joint approximation to any target $m-n$ and use the two Lipschitz inequalities through the intermediate witness. Their sum tends to zero.
\end{proof}

\begin{corollary}[Exact search at every precision]\label{cor:JA-search}
For every $\delta>0$ and every fixed integer target $t$, there are arbitrarily large positive integers $k$ with
\[
P_3(k)<\delta,\qquad P_2(k-t)<\delta.
\]
Consequently the exhaustive certified interval search of Section~\ref{sec:interval-search} terminates at every prescribed strict precision and target.
\end{corollary}
\begin{proof}
Obtain one witness at precision $\delta/2$ from Proposition~\ref{prop:JA-equivalent}. By Lemma~\ref{lem:nondiscrete}, there are positive mixed-null integers tending to infinity: change signs of nonzero null points, and note that a fixed finite set cannot be null. Add sufficiently large such integers to the witness and use subadditivity. The search uses exact finite comparisons and full-tail bounds; a positive witness therefore eventually lies in its enumerated range.
\end{proof}
This is an existence and termination result for simultaneous modular constraints at sparse square positions. It does not provide an effective practical upper bound for the least witness, polynomial runtime, a limiting distribution, or independence of base-$2$ and base-$3$ digits. The proof uses smoothing, finite Fourier orthogonality, normalized convolution, and a Diophantine separation estimate. The unweighted collision energies of the original centred balls need not satisfy the same uniform bound as the specially translated smooth weights. Optimizing the exact search and certifying large floating spectral experiments remain separate numerical tasks.

\begin{corollary}[Every fixed prime pair]\label{rem:JA-other-pairs}
For every fixed pair of distinct primes $p,q$, $\mathrm{JA}^{p,q}_t$ holds for every integer target $t$, and
\[
H_{\{p,q\}}=\Z[1/(pq)]/\Z.
\]
The completion of the two-prime metric is $G_p\times G_q$. In particular no irrational parameter is continuous for any such pair.
\end{corollary}
\begin{proof}
Replace $2,3$ by $p,q$ in Sections~\ref{sec:normalized-product} and \ref{sec:random-JA}. The individual window bounds and finite Fourier identity do not depend on the primes. Their successive square-level ratios still grow geometrically. The merged growth argument holds because $2\sqrt{\log p}$ and $2\sqrt{\log q}$ are both greater than one. Coprimality excludes exact two-term cancellation, and Matveev's estimate applies to the two fixed prime logarithms and the small rational coefficient ratio, with constants depending on $p,q$. The random-translation, support and packing arguments are unchanged; the quotient and bounded-gap proofs use only the corresponding one-prime duals. This proves $\mathrm{JA}^{p,q}_1$. Apply Theorem~\ref{thm:general-pair} for all integer targets and the dual, and the diagonal-density argument of Proposition~\ref{prop:JA-equivalent} for the completion.
\end{proof}
For three or more primes, simultaneous approximation to independent targets in every prime factor remains to be established. A direct product-quotient argument or an approximation theorem for multi-prime blocks could address it; pairwise joint approximation alone does not give the finite-family conclusion.

\begin{corollary}[Compatible joins of real TVG topologies]\label{cor:real-TVG-join}
Let $f:\R\to\Q$ be a nonzero $\Q$-linear map. For a prime $p$, let $\gamma_p$ be the group topology on $\Q$ generated by
\[
P_{p,d}(u)=\sup_{n\ge0}\nrm{u/(dp^{n^2})},
\qquad d\in\Z,\ d\ge1.
\]
Let $\eta_p=T_f(\gamma_p)$ be generated by the usual real coordinate and all maps $f_a(x)=f(ax)$ into $(\Q,\gamma_p)$, with $a\in\R$. Here a TVG topology means a group topology for which every fixed real scalar multiplication is continuous. For distinct primes $p,q$, these are Hausdorff locally quasi-convex TVG topologies and
\[
\widehat{(\Q,\gamma_p\vee\gamma_q)}=\Gamma,
\qquad \Gamma=\{u\mapsto ru+\Z:r\in\Q\},
\]
\[
\eta_p\vee\eta_q=T_f(\gamma_p\vee\gamma_q),
\qquad \widehat{(\R,\eta_p\vee\eta_q)}=H_f,
\]
\[
H_f=\{x\mapsto cx+f(ax)+\Z:c,a\in\R\}.
\]
Each $\eta_p$ has this same full character group $H_f$. The join is strictly finer than each of $\eta_p$ and $\eta_q$.
\end{corollary}
\begin{proof}
The one-prime compatibility and the rational-transfer theorem of \cite{Previous} apply to these topologies. The integer-polar reduction there uses only the $n=0$ observation and a common denominator, so applies with any distinct primes $p,q$. Corollary~\ref{rem:JA-other-pairs} says that the integer dual is $\Z[1/(pq)]/\Z$, which consists of rational parameters. The reduction therefore gives the displayed rational-join dual $\Gamma$. The transfer theorem preserves joins and recovers the rational topology as the open quotient under $f$. It gives the displayed real-join identity, the full dual $H_f$, and the stated topological properties. All real indices $a$ in the definition of $T_f$ are used here.

For strictness, it suffices by order reflection of the transfer to check that $\gamma_p$ and $\gamma_q$ are incomparable. Fix $d\ge1$, write $d=p^v e$ with $p\nmid e$, choose $n^2>v$, and set
\[
m=p^{n^2-1-v},\qquad u_k=dmq^{k^2}.
\]
Then
\[
P_{q,d}(u_k)\le m q^{-2k-1}\longrightarrow0,
\qquad \nrm{u_k/p^{n^2}}
=\nrm{eq^{k^2}/p}\ge1/p.
\]
Any finite collection of $q$-prime observations is controlled by one such denominator $d$, by taking a common multiple. Thus no $\gamma_q$-neighbourhood is contained in the fixed $\gamma_p$-neighbourhood $P_{p,1}<1/(2p)$. Interchanging the primes proves incomparability. The same holds for their transfers, so their join strictly contains each.
\end{proof}

\begin{remark}[Integer parameters and real-line characters]\label{rem:integer-real-parameters}
An irrational parameter on the integer group means a character $k\mapsto k\alpha+\Z$ with $\alpha\notin\Q$; the integer $k$ is its argument. Corollary~\ref{rem:JA-other-pairs} excludes such parameters for the two-prime integer metric. It does not exclude irrational coefficients of ordinary real-line characters. For example
\[
x\mapsto\sqrt{2}x+\Z
\]
is continuous in each $\eta_p$ and in their join. On the other hand
\[
x\mapsto\sqrt{2}f(x)+\Z
\]
is not continuous in that join: since $f$ is an open quotient onto $(\Q,\gamma_p\vee\gamma_q)$, its continuity would make $u\mapsto\sqrt{2}u+\Z$ continuous there, contradicting the rational dual $\Gamma$. The character $x\mapsto f(\sqrt{2}x)+\Z$ is continuous, since its irrational input rescaling is one of the defining observations.
\end{remark}
These are specific compatible joins with a common prescribed real dual. They do not alter the incompatible two-chain example or the resulting nonexistence of a greatest compatible topology in Part~I of \cite{HamelPair}. They do not compute nuclearity of the real transfers. The product-completion formula in Corollary~\ref{rem:JA-other-pairs} concerns the auxiliary integer metric; no product-completion formula for the full real transferred topology is asserted here.

\begin{proposition}[Character extension criterion]\label{prop:mixed-extension}
Let $p,q$ be distinct primes and let $E=c_0(\T)$ carry the uniform metric $d_E(u,v)=\sup_{n\ge1}\nrm{u_n-v_n}$. Define
\[
J(k)=\bigl((k/p^{n^2}+\Z)_{n\ge1},\ (k/q^{n^2}+\Z)_{n\ge1}\bigr).
\]
Then $J$ is an isometric group embedding of $(\Z,P_{\{p,q\}})$ into $E\times E$, with its maximum product metric. Moreover,
\[
H_{\{p,q\}}=\Z[1/(pq)]/\Z
\]
if and only if $J(\Z)$ is dually embedded in $E\times E$, meaning that every continuous character of $J(\Z)$ extends to a continuous character of $E\times E$.
\end{proposition}
\begin{proof}
For each fixed $k$, both coordinate sequences tend to zero. Their uniform distances give exactly $P_{\{p,q\}}(k-k')$, and this metric separates integers, proving the embedding assertion. Continuous characters of $c_0(\T)$ are finite integer combinations of coordinate evaluations, as used in \cite[Corollary~4.5]{ABM}. The characters of $E\times E$ are therefore the sums of two such combinations. Their restrictions to $J(\Z)$ have exactly the parameters
\[
\Z[1/p]/\Z+\Z[1/q]/\Z=\Z[1/(pq)]/\Z.
\]
The last equality follows from B\'ezout's identity, and every term is realized by finitely many square coordinates. If $J(\Z)$ is dually embedded, every continuous mixed character has one of these parameters. Conversely, if the displayed mixed-dual equality holds, each continuous character has such a finite-coordinate representation and hence extends. This proves the equivalence.
\end{proof}
This criterion does not assume joint approximation. Theorem~\ref{thm:joint-approximation} and Corollary~\ref{rem:JA-other-pairs} establish its dual equality in the present setting, and hence the extension property. It is a concrete reformulation of the arithmetic conclusion, rather than an independent proof of it. The related general connection between compatible joins and diagonal character extension is \cite[Proposition~3.9]{GabCompatible}; that proposition concerns topologies for one fixed duality. Here the separate integer prime-square topologies have different duals, so the preceding proof uses the coordinate embedding directly.

For comparison, let $b_n=p^n$. The topology $\Gamma_b$ in \cite[Definition~3.1]{DBM2017} is the supremum over uniform-convergence topologies associated with sparse subsequences of this chain. Its dual is $\Z[1/p]/\Z$ by Lemma~3.5 there, while Theorem~3.3 shows that it has no nontrivial convergent sequences and is not metrizable. Thus a compatible supremum in that prime-adic setting is already known. It neither supplies a greatest compatible topology here nor replaces the mixed-prime approximation argument.

\section{Remaining questions}\label{sec:open}
Joint approximation and the full continuous character group are now determined for every fixed pair of distinct primes. The following quantitative and finite-family questions remain.
\begin{enumerate}
\item Obtain useful bounds for the least positive strict witness as a function of the precision and target. The existence proof supplies no practical complexity estimate for the certified search. The earlier late-survival hypothesis is still a separate quantitative statement; joint approximation does not prove it.
\item Determine whether the original unweighted, zero-centred collision energies admit uniform upper bounds. The translated smooth-weight construction is sufficient for joint approximation and does not establish this stronger property of the original balls.
\item Determine the exact neighbourhood gcds and fine-scale counting laws. Restricting their possible prime divisors does not compute their valuations, densities or correlations.
\item For three or more primes, determine the full mixed dual and prove or disprove simultaneous approximation to independent targets in all prime factors. A direct approach through the product quotient need not proceed by partitions into blocks. Pairwise approximation alone is not asserted to imply the finite-family statement or the sum-of-pairs description in Question~\ref{q:additive}.
\item Remove the remaining hypothesis of the Lean formalization of joint approximation, by formalizing Matveev's bound for linear forms in logarithms and its conversion to the two-term separation estimate \eqref{eq:two-term-separation}. The rest of the argument, from the window kernels to joint approximation, is checked in Lean conditionally on that estimate; the scope is described in the section on checks by Lean and Mathlib.
\end{enumerate}

\section*{Checks by Lean and Mathlib}
\phantomsection
\addcontentsline{toc}{section}{Checks by Lean and Mathlib}
\emph{Formalization of joint approximation.} Theorem~\ref{thm:joint-approximation} has been formalized in Lean~4 with Mathlib. The main formal statement, \texttt{JA\_of\_matveev}, says that for every $A\ge0$ the two-term separation estimate \eqref{eq:two-term-separation} implies joint approximation, in the form that for every $\delta>0$ there is an integer $k$ with $P_3(k)<\delta$ and $P_2(k-1)<\delta$. The estimate is stated in Lean as the hypothesis \texttt{TwoTermSeparation A}, for pairs of prime-square periods of $2$ and $3$ that are close in the sense $Q<Q'<tQ$. All other steps of the argument are checked: the window kernel and its tail bound, the positive cascade obtained by random relative translations with its mass bounds, the product identity for the normalized factors and the center maximality, the far-field recursion and the growth of merged square powers, the estimates of the factors at separated, close and non-deep levels (Section~\ref{sec:normalized-product}), the final product bound and the assembly of the central bound, the energy and packing arguments, the transfer to joint approximation (Section~\ref{sec:random-JA}), and the passage from syndeticity of all $W_\varepsilon$ to joint approximation (Section~\ref{sec:quotient}).

The last step is checked by an argument different from the printed proof of Theorem~\ref{thm:quotient-obstruction}. The one-prime character theorem, Theorem~\ref{thm:oneprime}, is proved in Lean. The existence of a nontrivial continuous character on a precompact seminormed copy of $\Z$ is obtained from the spectral theorem for a compact self-adjoint convolution operator on $L^2$ of the Haar measure of the completion, instead of Pontryagin duality for the quotient group. The full statements of Theorem~\ref{thm:quotient-obstruction}, namely the Polish group structure and the equivalence of the compactness conditions, are not formalized.

All declarations of this development depend only on the standard axioms \texttt{propext}, \texttt{Classical.choice} and \texttt{Quot.sound}, and the files contain no \texttt{sorry}. The scope is as follows. The hypothesis \texttt{TwoTermSeparation} is not proved. Matveev's bound and its conversion to \eqref{eq:two-term-separation} are not formalized, so joint approximation is formalized conditionally and not as an unconditional Lean theorem. The formal statements follow the structure of the printed proof, but some steps are organized differently; for example, the central bound is assembled for the expected cascade measures at buffered heights. The correspondence of each formal statement with the printed statement still requires review by a reader.

Selected statements have been checked using Lean and Mathlib. The compactness and Lipschitz splitting implication in Theorem~\ref{thm:JA-split} has been independently typechecked with the joint approximation limits as hypotheses; it does not prove those limits or joint approximation. Earlier checked statements include: the Lipschitz implication and flatness equivalence, valuation bounds and the stated exact radii, the coloured-difference equivalence, profile summation and the finite-family bound with the two discontinuity tests, pointwise sparse-tail estimates and the one-sided exclusion criterion, and the finite-scale upper counting estimate used in Theorem~\ref{thm:density}. These checks concern the formal statements under their stated hypotheses, rather than every printed proof. They do not cover the density lower bound or its full asymptotic specialization, the stationary-inner construction, the thickness argument and cited transversality results, the limiting step in Proposition~\ref{prop:sparse-tests}(a), or the near-resonance example. The new nonzero-target reduction, residue-distribution applications, block estimates, and survival-count arguments have not been checked in Lean. The finite survivor table was independently reproduced using exact integer comparisons. The prime-pair and finite-prime generalizations in Section~\ref{sec:general-primes}, the interval construction, and its numerical certificates have not been checked in Lean. The computations in Section~\ref{sec:interval-search} use exact integer or rational arithmetic, independent replay of interval certificates, an independent implementation without a CRT prefix, and a repeated enumeration in both orientations. These are computational checks, not formal Lean proofs. Of the results in Section~\ref{sec:quotient}, only the implication in Corollary~\ref{cor:syndetic-JA} from syndeticity to joint approximation is checked in Lean, by the argument described above; the remaining quotient, compactness, covering and finite-coordinate statements are not. Its shifted-target witnesses were independently verified with full-tail bounds, without a formal Lean proof. The density, recursion, buffered-height and energy arguments in Sections~\ref{sec:density-JA}--\ref{sec:buffered-energy} are used, in the form needed there, in the formal proof of joint approximation described above; statements in those sections that are not used there, including the common-Lipschitz criteria, have not been checked in Lean. Exact pair-counting tests included direct enumeration at smaller heights and a separate implementation of the large-height certificate. The quadratic state count and the periodic-kernel deductions are mathematical arguments; these computational checks do not formally verify them or establish any uniform upper estimate.

The Lean and Mathlib sources for the selected formalized results are available in the public repository \href{https://github.com/irinayar/MyPseudotopology-Lean}{MyPseudotopology-Lean}. The source snapshot cited here is \href{https://github.com/irinayar/MyPseudotopology-Lean/tree/c0444e290d5b7f0952907356f7f6a23e62b87586}{commit \texttt{c0444e2}}. Its verification records state the checked declarations and their source hashes. Individual module checks are distinguished from a full rebuild of the entire library. The formalization of joint approximation is in \href{https://github.com/irinayar/MyPseudotopology-Lean/tree/2361e26997526a7b0dcdc4c221ed3768955c96ab}{commit \texttt{2361e26}}; its file \texttt{JA\_FORMALIZATION.md} gives the module map and \texttt{verification/ja-check.txt} lists the axioms of the main declarations. The 47 modules \texttt{JA*.lean} were built together with \texttt{lake build}; the root library target \texttt{LeanCheck} was not rebuilt for this commit.

Checked using Lean and Mathlib: the abstract far-field recursion and its geometric specialization, finite cyclic Fourier autocorrelation and positivity, finite cosine-sum center maximality, product telescoping, and coprime two-term nonvanishing. These seven earlier declarations concern the stated abstract lemmas under their explicit hypotheses.

Thirteen additional declarations have been independently typechecked. They establish absolute summability and center maximality for nonnegative-coefficient cosine series; a finite-product bound with boundedly many exceptional factors; lower and two-sided product bounds for relative normalization errors; geometric sums and their paired finite-prefix bounds; a numerical square-level spacing inequality and the two-colour three-level span argument; the finite packing-to-cover principle for a symmetric irreflexive relation with uniformly bounded finite packings; and exact finite cyclic-character bilinear orthogonality for translation averaging. All thirteen declarations use only the standard axioms of propositional extensionality, classical choice and quotient soundness, with no admitted proofs or added axioms.

The thirteen declarations above were checked before the formalization of joint approximation described at the beginning of this section. The later formalization reuses several of them, for example the far-field recursion, the finite-product bound with exceptional factors, the three-level span argument and the packing-to-cover principle, and adds the window and kernel construction, the tail estimate, the deep-level and terminal-level estimates, the spatial lift-sum bounds, the buffered-height and spectral-energy transfer, and the character-separation argument. What is not formalized is the quantitative Matveev separation and its conversion to \eqref{eq:two-term-separation}, which enters the formal statement as a hypothesis. The splitting implication in Theorem~\ref{thm:JA-split} remains a separate statement that keeps joint approximation as a hypothesis.

Finite FFT diagnostics are computational checks, not interval certificates or formal proofs. The finite weighted-count comparison in Lemma~\ref{lem:finite-weighted-count} and the real-TVG consequences in Corollary~\ref{cor:real-TVG-join} and Remark~\ref{rem:integer-real-parameters} have been checked mathematically, but have not been checked in Lean.

The character-extension reformulation in Proposition~\ref{prop:mixed-extension} and the applications of the cited quasi-convexity and sparse-supremum results have not been checked in Lean. They add no new formalization claim.

\section*{Acknowledgement}
\phantomsection
\addcontentsline{toc}{section}{Acknowledgement}
The author acknowledges OpenAI's Codex and Anthropic's Claude for assistance in mathematical discussion, review of arguments, exposition, and article preparation. The author also gratefully acknowledges the developers and contributors of the Lean proof assistant and the mathlib library, which were used for formal verification of selected results and of the conditional joint approximation theorem. The author assumes responsibility for the final content and any remaining errors.

\end{document}